\documentclass[11pt, a4paper]{article}
\usepackage{latexsym,amssymb}
\usepackage{amsmath}
\usepackage[mathscr]{eucal}
\usepackage{color}
\usepackage[abbrev]{amsrefs}
\usepackage{ulem}

\newtheorem{Theorem}{\bf Theorem}[section]
\newtheorem{Lemma}{\bf Lemma}[section]
\newtheorem{Proposition}{\bf Proposition}[section]
\newtheorem{Corollary}{\bf Corollary}[section]
\newtheorem{Remark}{\bf Remark}[section]
\newtheorem{Example}{\bf Example}[section]
\newtheorem{Definition}{\bf Definition}[section]

\newenvironment{theorem}{\begin{Theorem}$\!\!\!$}{\end{Theorem}}
\newenvironment{lemma}{\begin{Lemma}$\!\!\!$}{\end{Lemma}}
\newenvironment{proposition}{\begin{Proposition}$\!\!\!$}{\end{Proposition}}
\newenvironment{corollary}{\begin{Corollary}$\!\!\!$}{\end{Corollary}}
\newenvironment{remark}{\begin{Remark}$\!\!\!$}{\end{Remark}}

\newenvironment{definition}{\begin{Definition}$\!\!\!$}{\end{Definition}}

\numberwithin{equation}{section}

\numberwithin{equation}{section}

\begin{document}
\title{Existence of solutions to the fractional semilinear heat\\ equation 
with a singular inhomogeneous term}
\author{Kazuhiro Ishige, Tatsuki Kawakami, and Ryo Takada}
\date{}
\maketitle
\begin{abstract}
We study the existence of solutions to 
the fractional semilinear heat equation with a singular inhomogeneous term. 
For this aim, 
we establish decay estimates of the fractional heat semigroup in several uniformly local Zygumnd spaces. 
Furthermore, we apply the real interpolation method in uniformly local Zygmund spaces to obtain sharp integral estimates 
on the inhomogeneous term and the nonlinear term. 
This enables us to find sharp sufficient conditions 
for the existence of solutions to the fractional semilinear heat equation with a singular inhomogeneous term.
\end{abstract}
%

\vspace{40pt}
\noindent 
Addresses:

\smallskip
\noindent 
K. I.: Graduate School of Mathematical Sciences, The University of Tokyo,\\ 
3-8-1 Komaba, Meguro-ku, Tokyo 153-8914, Japan\\
\noindent 
E-mail: {\tt ishige@ms.u-tokyo.ac.jp}\\

\smallskip
\noindent 
{T. K.}: Applied Mathematics and Informatics Course,\\ 
Faculty of Advanced Science and Technology, Ryukoku University,\\
1-5 Yokotani, Seta Oe-cho, Otsu, Shiga 520-2194, Japan\\
\noindent 
E-mail: {\tt kawakami@math.ryukoku.ac.jp}\\

\smallskip
\noindent 
R. T.: Graduate School of Mathematical Sciences, The University of Tokyo,\\ 
3-8-1 Komaba, Meguro-ku, Tokyo 153-8914, Japan\\
\noindent 
E-mail: {\tt r-takada@g.ecc.u-tokyo.ac.jp}\\

\newpage
\section{Introduction}
Consider the Cauchy problem for the fractional semilinear heat equation with an inhomogeneous term
\begin{equation}
\tag{P}
\label{eq:P}
	\left\{
	\begin{array}{ll}
	\partial_t u+(-\Delta)^{\frac{\theta}{2}}u=|u|^{p-1}u+\mu, & x\in{\mathbb R}^N,\,\,t>0,\vspace{5pt}\\
	u(x,0)=0, & x\in{\mathbb R}^N,
	\end{array}
	\right.
\end{equation}
where $N\ge 1$, $\partial_t:=\partial/\partial t$, $\theta\in(0,2]$, $p>1$, 
and $\mu$ is a (possibly signed) Radon measure or a locally integrable function in ${\mathbb R}^N$. 
Here $(-\Delta)^{\theta/2}$ denotes the fractional power of the Laplace operator $-\Delta$ in ${\mathbb R}^N$. 
Problem~\eqref{eq:P} has been studied in many papers (see e.g., \cites{BLZ, CV, CVW, HIT, KK01, KK02, KQ, KT, Zeng, Zhang01, Zhang02, Zhang03} and a survey book~\cite{QS}) 
and appears in the study of the super-Brownian motion (see \cite{Lee}). 
In this paper we obtain sharp sufficient conditions on the inhomogeneous term $\mu$ 
for the existence of solutions to problem~\eqref{eq:P}.

This paper is motivated by~\cite{HIT}, which gave necessary conditions and sufficient conditions 
on the inhomogeneous term $\mu$ for the existence of nonnegative solutions to problem~\eqref{eq:P} 
under the following mild definition of solutions. 
\begin{definition}
\label{Definition:1.1}
Let $N\ge 1$, $\theta\in(0,2]$, $p>1$, $T\in(0,\infty]$, and $\mu$ be a Radon measure in ${\mathbb R}^N$. 
For any measurable function $u$ in ${\mathbb R}^N\times(0,T)$, 
we say that $u$ is a solution to problem~\eqref{eq:P} in ${\mathbb R}^N\times[0,T)$ if  
$|u|^{p-1}u$ is locally integrable in ${\mathbb R}^N\times[0,T)$ and $u$ satisfies 
\begin{equation*}
\int_0^T\int_{{\mathbb R}^N} u\left(-\partial_t \varphi + (-\Delta)^\frac{\theta}{2} \varphi\right)\,dx\,dt = 
\int_0^T\int_{{\mathbb R}^N}  |u|^{p-1}u\varphi\,dx\,dt + \int_0^T \int_{{\mathbb R}^N}  \varphi \,d\mu(x)\,dt
\end{equation*}
for $\varphi\in C^\infty_0({\mathbb R}^N\times[0,T))$.
\end{definition}
Under Definition~\ref{Definition:1.1},  
Hisa~{\it et al.}~\cite{HIT} developed the arguments in \cites{HI01, KK01, KK02} 
to obtain necessary conditions on the inhomogeneous term $\mu$ 
for the existence of nonnegative solutions to problem~\eqref{eq:P} (see \cite{HIT}*{Theorem~1.1}).
Set $B(x,\rho):=\{y\in\mathbb R^N : |x-y|<\rho\}$ for $x\in\mathbb R^N$ and $\rho\in(0,\infty)$.
\begin{theorem}
\label{Theorem:1.1}
Let $N\ge 1$, $\theta\in(0,2]$, and $p>1$. 
Let $\mu$ be a nonnegative Radon measure in ${\mathbb R}^N$. 
Assume that problem~\eqref{eq:P} possesses a nonnegative solution in ${\mathbb R}^N\times[0,T)$ 
for some $T\in(0,\infty)$. 
Set 
$$
p_*:=\frac{N}{N-\theta}\quad\mbox{if}\quad N>\theta,
\qquad
p_*:=\infty\quad\mbox{if}\quad N\le\theta. 
$$
\begin{itemize}
  \item[{\rm (1)}] {\rm ({\it Subcritical Case})}
  If $1<p<p_*$, then there exists $\gamma_1>0$ such that 
  $$
  \sup_{z\in{\mathbb R}^N}\,\mu(B(z,T^{\frac{1}{\theta}}))\le\gamma_1 T^{\frac{N}{\theta}-\frac{p}{p-1}}.
  $$
  \item[{\rm (2)}] {\rm ({\it Critical Case})}  
  If $p=p_*$, then there exists $\gamma_2>0$ such that 
  $$
  \sup_{z\in{\mathbb R}^N}\,\mu(B(z,\rho))\le\gamma_2\left[\log\left(e+\frac{T^{\frac{1}{\theta}}}{\rho}\right)\right]^{-\frac{N}{\theta}+1},
  \quad \rho\in(0,T^{\frac{1}{\theta}}).
  $$
  \item[{\rm (3)}] {\rm ({\it Supercritical Case})}  
  If $p>p_*$, then there exists $\gamma_3>0$ such that 
  $$
  \sup_{z\in{\mathbb R}^N}\,\mu(B(z,\rho))\le\gamma_3\rho^{N-\frac{\theta p}{p-1}},
  \quad \rho\in(0,T^{\frac{1}{\theta}}).
  $$
\end{itemize} 
\end{theorem}
\begin{remark}
\label{Remark:1.1}
Let $1<p\le p_*$ and $\mu$ be a nonnegative Radon measure in ${\mathbb R}^N$. 
Assume that problem~\eqref{eq:P} possesses a global-in-time nonnegative solution. 
Then, by assertions~{\rm (1)} and {\rm (2)} we see that $\mu=0$ in ${\mathbb R}^N$. 
Indeed, if $1<p<p_*$, then it follows from assertion~{\rm (1)} that 
$$
0\le\mu({\mathbb R}^N)=\lim_{T\to\infty}\sup_{z\in{\mathbb R}^N}\,\mu(B(z,T^{\frac{1}{\theta}}))
\le\gamma_1\lim_{T\to\infty}T^{\frac{N}{\theta}-\frac{p}{p-1}}=0.
$$
Similarly, if $p=p_*$, then it follows from assertion~{\rm (2)} that 
\begin{equation}
\label{eq:1.1}
0\le\mu({\mathbb R}^N)=\lim_{T\to\infty}\sup_{z\in{\mathbb R}^N}\,\mu(B(z,T^{\frac{1}{2\theta}}))
\le\gamma_2\lim_{T\to\infty}\left[\log\left(e+T^{\frac{1}{2\theta}}\right)\right]^{-\frac{N}{\theta}+1}=0.
\end{equation}
\end{remark}
On the other hand, in \cite{HIT}
sufficient conditions on the inhomogeneous term $\mu$ 
for the existence of nonnegative solutions to problem~\eqref{eq:P}
were proved by the so-called supersolution methods (see \cite{HIT}*{Theorems~1.2 and 1.4}). 
Then, employing the arguments in \cite{TW} (see also \cites{IKO01, IKO02}),
we easily extend the results in \cite{HIT}*{Theorems~1.2 and 1.4} to 
obtain the following theorem.
\begin{theorem}
\label{Theorem:1.2}
Let $N\ge 1$, $\theta\in(0,2]$, and $p>1$. 
Then there exists $\varepsilon>0$ with the following properties. 
\begin{itemize}
  \item[{\rm (4)}] {\rm ({\it Subcritical Case})}
  If $1<p<p_*$ and a Radon measure $\mu$ in ${\mathbb R}^N$ satisfies
  $$
  \sup_{z\in{\mathbb R}^N}\,|\mu|(B(z,T^{\frac{1}{\theta}}))\le\varepsilon T^{\frac{N}{\theta}-\frac{p}{p-1}}
  $$
  for some $T\in(0,\infty)$, then problem~\eqref{eq:P} possesses a solution in ${\mathbb R}^N\times[0,T)$. 
  Here $|\mu|(B(z,\rho))$ is the total variation of $\mu$ in $B(z,\rho)$, where $z\in{\mathbb R}^N$ and $\rho\in(0,\infty)$.
  \item[{\rm (5)}] {\rm ({\it Critical Case})}  
  If $p=p_*$ and a measurable function $\mu$ in ${\mathbb R}^N$ satisfies
  $$
  |\mu(x)|\le \varepsilon |x|^{-N}\left[\log\left(e+\frac{1}{|x|}\right)\right]^{-\frac{N}{\theta}}+C_1,\quad \mbox{a.a.~$x\in{\mathbb R}^N$},
  $$
  for some $C_1\ge 0$, then problem~\eqref{eq:P} possesses a local-in-time solution.
  \item[{\rm (6)}] {\rm ({\it Supercritical Case})}  
  If $p>p_*$ and a measurable function $\mu$ in ${\mathbb R}^N$ satisfies
  $$
  |\mu(x)|\le\varepsilon |x|^{-\frac{\theta p}{p-1}}+C_2,\quad \mbox{a.a.~$x\in{\mathbb R}^N$},
  $$
  for some $C_2>0$ {\rm ({\it resp. $C_2=0$})}, then problem~\eqref{eq:P} possesses a local-in-time solution 
  {\rm ({\it resp. global-in-time solution})}.
\end{itemize} 
\end{theorem}
In the subcritical case, 
by assertions~(1) and (4)
we see that, for any nonnegative Radon measure~$\mu$ in ${\mathbb R}^N$, 
problem~\eqref{eq:P} possesses a local-in-time nonnegative solution if and only if 
$$
\sup_{x\in{\mathbb R}^N}\mu(B(x,1))<\infty.
$$
In the critical and the supercritical cases, 
Theorem~\ref{Theorem:1.1}~(2), (3) and Theorem~\ref{Theorem:1.2}~(5), (6) imply that 
the strength of the singularity at the origin of the function $\mu_c$ in ${\mathbb R}^N$ defined by 
\begin{equation}
\label{eq:1.2}
\mu_c(x):=
\left\{\begin{array}{ll}
 \displaystyle{|x|^{-N}\left[\log\left(e+\frac{1}{|x|}\right)\right]^{-\frac{N}{\theta}}}, & \qquad\mbox{if}\qquad p=p_*,\vspace{7pt}\\
 |x|^{-\frac{\theta p}{p-1}}, & \qquad\mbox{if}\qquad p>p_*,\vspace{3pt}
 \end{array}\right.
\end{equation}
is a threshold for the local solvability of problem~\eqref{eq:P} in the following senses.
\begin{itemize}
  \item 
  If a nonnegative measurable function $\mu$ in ${\mathbb R}^N$ satisfies $\mu_c=o(\mu)$ uniformly near the origin, 
  then problem~\eqref{eq:P} possesses no local-in-time nonnegative solutions; 
  conversely, if $\mu =o(\mu_c)$ uniformly near the origin and $\mu$ is bounded away from the origin, 
  then problem~\eqref{eq:P} possesses a local-in-time nonnegative solution. 
  \item 
  Let $\mu=\lambda\mu_c$ with $\lambda>0$. 
  Then problem~\eqref{eq:P} possesses a local-in-time nonnegative solution if $\lambda$ is small enough; 
  problem~\eqref{eq:P} possesses no local-in-time nonnegative solutions if $\lambda$ is large enough. 
  In fact, one can easily show via supersolution methods (see \cite{HIT}*{Lemma~4.1}) 
  that there exists a unique $\lambda_*>0$ such that 
  problem~\eqref{eq:P} with $\mu =\lambda \mu_c$ possesses a local-in-time nonnegative solutions for all $\lambda\in (0,\lambda_*)$ 
  while it possesses no local-in-time nonnegative solutions for any $\lambda>\lambda_*$. 
\end{itemize}
The strength of the singularity of $\mu_c$ is useful to find sharp sufficient conditions 
for the existence of solutions to problem~\eqref{eq:P}. 
On the other hand, 
the argument in the proof of \cite{HIT}*{Theorem~1.2}
heavily depends on the explicit expression of the function $\mu_c$, 
and they are not applicable to the study of sharp sufficient conditions on the existence of solutions to problem~\eqref{eq:P} 
with more general singular inhomogeneous terms.  

In this paper we find a Banach space $X$ with $\mu_c\in X$ such that 
problem~\eqref{eq:P} is locally solvable for small $\mu\in X$ in the critical and the supercritical cases, respectively, 
using uniformly local Zygmund spaces. We establish decay estimates of the fractional heat semigroup in several uniformly local Zygumnd spaces. 
Then we obtain sharp sufficient conditions on the existence of solutions to problem~\eqref{eq:P} in the critical and the supercritical cases 
using the real interpolation method in uniformly local Zygumnd spaces. 
\subsection{Zygmund spaces}
We introduce some notations. 
We denote by ${\mathcal M}$ the set of Lebesgue measurable sets in ${\mathbb R}^N$. 
For any $E\in{\mathcal M}$, we denote by $|E|$ and $\chi_E$ the $N$-dimensional Lebesgue measure of $E$ 
and the characteristic function of $E$, respectively. 
Let ${\mathcal L}$ be the set of measurable functions in ${\mathbb R}^N$.
For any $q\in[1,\infty]$, we denote by $L^q$ and $\|\cdot\|_{L^q}$ the usual $L^q$-space on ${\mathbb R}^N$ and its norm, respectively.
For any $f\in{\mathcal L}$, 
we denote by $d_f$ the distribution function of~$f$, that is, 
$$
d_f(\lambda):=\left|\{x\in\mathbb R^N\,:\,|\,f(x)|>\lambda\}\right|, 
\quad \lambda > 0.
$$
We define the non-increasing rearrangement $f^*$ of $f$ by 
$$
f^{*}(s):=\inf\{\lambda>0\,:\,d_f(\lambda)\le s\},\quad s\in[0,\infty).
$$
Here we adopt the convention $\inf\emptyset=\infty$.
Then $f^*$ is right continuous in $[0,\infty)$.

Let $q\in[1,\infty]$ and $\alpha\in[0,\infty)$. Set 
\begin{equation}
\label{eq:1.3}
\Phi(s) :=\log(e+s),\quad s\in[0,\infty).
\end{equation}
The weak $L^q$-space $L^{q,\infty}$ is defined by 
$$
L^{q,\infty}:=\{f\in{\mathcal L}\,:\,\|f\|_{L^{q,\infty}}<\infty\},
\mbox{ where }
\|f\|_{L^{q,\infty}}:=
\left\{
	\begin{array}{ll}
	\displaystyle{\sup_{s>0}}\,\left\{sf^*(s)^q\right\}^{\frac{1}{q}}& \mbox{if}\,\,\, q\in[1,\infty),
	\vspace{5pt}\\
	\|f\|_{L^\infty} 
	& \mbox{if}\,\,\, q=\infty.\vspace{3pt}
	\end{array}
	\right.
$$
(See e.g., \cite{Grafakos}*{Section~1.4.2}.)
Define the Zygmund space $L^q(\log L)^\alpha$ by 
$$
	L^q(\log L)^\alpha :=\{f\in{\mathcal L}\,:\, \|f\|_{L^q(\log L)^\alpha}<\infty\},
$$
where
\begin{equation}
\label{eq:1.4}
	 \|f\|_{L^q(\log L)^\alpha}
	 :=
	 \left\{
	\begin{array}{ll}
	\displaystyle{\left(\int_0^\infty \Phi(s^{-1})^\alpha
	f^*(s)^q\,ds\right)^\frac{1}{q}}
	& \mbox{if}\,\,\, q\in[1,\infty),
	\vspace{7pt}
	\\
	\|f\|_{L^\infty} 
	& \mbox{if}\,\,\, q=\infty.\vspace{3pt}
	\end{array}
	\right.
\end{equation}
Then $L^q(\log L)^\alpha$ with $q>1$ is a Banach space (see Lemma~\ref{Lemma:3.4}). 
Furthermore, 
we introduce two weak Zygmund spaces $L^{q,\infty}(\log L)^\alpha$ and ${\mathfrak L}^{q,\infty}(\log {\mathfrak L})^\alpha$.
The space $L^{q,\infty}(\log L)^\alpha$ is the standard weak Zygmund space, which is defined by 
$$
	L^{q,\infty}(\log L)^\alpha :=\{f\in{\mathcal L}\,:\, \|f\|_{L^{q,\infty}(\log L)^\alpha}<\infty\},
$$
where 
\begin{equation}
\label{eq:1.5}
	\|f\|_{L^{q,\infty}(\log L)^\alpha}
	:=
	\left\{
	\begin{array}{ll}
	\displaystyle{\sup_{s>0}\,
	\bigg\{s\Phi(s^{-1})^\alpha f^*(s)^q\bigg\}^\frac{1}{q}}
	& \mbox{if}\,\,\, q\in[1,\infty),
	\vspace{7pt}
	\\
	\|f\|_{L^\infty} 
	& \mbox{if}\,\,\, q=\infty.\vspace{3pt}
	\end{array}
	\right.
\end{equation}
Then $L^{q,\infty}(\log L)^\alpha$ with $q>1$ is a Banach space (see Lemma~\ref{Lemma:3.4}). 
(See also e.g., \cite{BS}*{Chapter~$4$, Section~$6$} and \cite{Wadade}.)
The space ${\mathfrak L}^{q,\infty}(\log {\mathfrak L})^\alpha$ was introduced in \cite{IIK} and it is defined by 
$$
	{\mathfrak L}^{q,\infty}(\log {\mathfrak L})^\alpha 
	:=\{f\in {\mathcal L}\,:\,\|f\|_{{\mathfrak L}^{q,\infty}(\log{\mathfrak L})^\alpha}<\infty\},
$$
where 
\begin{equation}
\label{eq:1.6}
	\|f\|_{{\mathfrak L}^{q,\infty}(\log{\mathfrak L})^\alpha}
	:=
	\left\{
	\begin{array}{ll}
	\displaystyle{\sup_{s>0}\,\bigg\{\Phi(s^{-1})^\alpha
	\int_0^s f^*(\tau)^q\,d\tau\bigg\}^{\frac{1}{q}}}
	& \mbox{if}\,\,\, q\in[1,\infty),
	\vspace{7pt}
	\\
	\|f\|_{L^\infty} 
	& \mbox{if}\,\,\, q=\infty.\vspace{3pt}
	\end{array}
	\right.
\end{equation}
Then ${\mathfrak L}^{q,\infty}(\log{\mathfrak L})^\alpha$ is a Banach space equipped 
with the norm~$\|\,\cdot\,\|_{{\mathfrak L}^{q,\infty}(\log{\mathfrak L})^\alpha}$
(see \cite{IIK}*{Lemma~2.1}). 
Since $\Phi(s^{-1})$ and $f^*(s)$ are non-increasing for $s\in(0,\infty)$, we see that 
\begin{equation}
\label{eq:1.7}
\begin{split}
\|f\|_{L^q(\log L)^\alpha}
 & =\left(\int_0^\infty \Phi(s^{-1})^\alpha f^*(s)^q\,ds\right)^\frac{1}{q}\\
 & \ge \sup_{s>0}\,\bigg\{\Phi(s^{-1})^\alpha
 \int_0^s f^*(\tau)^q\,d\tau\bigg\}^{\frac{1}{q}}=\|f\|_{{\mathfrak L}^{q,\infty}(\log{\mathfrak L})^\alpha}\\
 & \ge\sup_{s>0}\,\bigg\{s\Phi(s^{-1})^\alpha
	f^*(s)^q\bigg\}^{\frac{1}{q}}=\|f\|_{L^{q,\infty}(\log L)^\alpha}
\end{split}
\end{equation}
for $f\in{\mathcal L}$, where $q\in[1,\infty)$ and $\alpha\in[0,\infty)$. This implies that 
$$
L^q(\log L)^\alpha\subset {\mathfrak L}^{q,\infty}(\log{\mathfrak L})^{\alpha}\subset L^{q,\infty}(\log L)^{\alpha},
$$
where $q\in[1,\infty]$ and $\alpha\in[0,\infty)$. 
Furthermore, 
\begin{equation}
\label{eq:1.8}
{\mathfrak L}^{q,\infty}(\log{\mathfrak L})^\alpha
\subsetneq
{\mathfrak L}^{q,\infty}(\log{\mathfrak L})^0=L^q(\log L)^0=L^q
\subsetneq L^{q,\infty}=L^{q,\infty}(\log L)^0,
\end{equation}
where $q\in[1,\infty)$ and $\alpha\in(0,\infty)$.

Next, we introduce uniformly local Zygmund spaces $L^q_{{\rm ul}}(\log L)^\alpha$, $L^{q,\infty}_{{\rm ul}}(\log L)^\alpha$,
and ${\mathfrak L}_{{\rm ul}}^{q,\infty}(\log{\mathfrak L})^\alpha$, where $q\in[1,\infty]$ and $\alpha\in[0,\infty)$, by
\begin{align*}
	L_{{\rm ul}}^q(\log L)^\alpha
	 & :=\{f\in{\mathcal L}\,:\,\|f\|_{L_{{\rm {ul}}}^q(\log L)^\alpha}<\infty\},\\
	L_{{\rm ul}}^{q,\infty}(\log L)^\alpha
	 & :=\{f\in{\mathcal L}\,:\,\|f\|_{L_{{\rm {ul}}}^{q,\infty}(\log L)^\alpha}<\infty\},
	\\
	{\mathfrak L}_{{\rm ul}}^{q,\infty}(\log{\mathfrak L})^\alpha
	 & :=\{f\in {\mathcal L}\,:\,\|f\|_{\mathfrak L_{{\rm {ul}}}^{q,\infty}(\log{\mathfrak L})^\alpha}<\infty\},
\end{align*}
where
\begin{align*}
	\|f\|_{L_{{\rm ul}}^q(\log L)^\alpha}
	 & :=\sup_{z\in{\mathbb R}^N}\,\|f\chi_{B(z,1)}\|_{L^q(\log L)^\alpha},\\
	\|f\|_{L_{{\rm ul}}^{q,\infty}(\log L)^\alpha}
	 & :=\sup_{z\in{\mathbb R}^N}\,\|f\chi_{B(z,1)}\|_{L^{q,\infty}(\log L)^\alpha},\\
	\|f\|_{{\mathfrak L}_{{\rm ul}}^{q,\infty}(\log{\mathfrak L})^\alpha}
	 & :=\sup_{z\in{\mathbb R}^N}\,\|f\chi_{B(z,1)}\|_{{\mathfrak L}^{q,\infty}(\log{\mathfrak L})^\alpha},
\end{align*}
respectively. 
We write 
$$
L_{{\rm ul}}^{q,\infty}:=L_{{\rm ul}}^{q,\infty}(\log L)^0, 
$$ 
which is a uniformly local weak Lebesgue space.  
For $f\in{\mathcal L}$ and $\rho\in(0,\infty]$, we often write
\begin{equation}
\label{eq:1.9}
\begin{split}
 & |f|_{q,\alpha;\rho}:=\sup_{z\in{\mathbb R}^N}\,\|f\chi_{B(z,\rho)}\|_{L^q(\log  L)^\alpha},\\
 & \|f\|_{q,\alpha;\rho}:=\sup_{z\in{\mathbb R}^N}\,\|f\chi_{B(z,\rho)}\|_{L^{q,\infty}(\log  L)^\alpha},
 \quad |||f|||_{q,\alpha;\rho}
:=\sup_{z\in{\mathbb R}^N}\,\|f\chi_{B(z,\rho)}\|_{{\mathfrak L}^{q,\infty}(\log{\mathfrak L})^\alpha}.
\end{split}
\end{equation}
\subsection{Main results}
We state our main results, Theorems~\ref{Theorem:1.3} and \ref{Theorem:1.4}. 
Theorem~\ref{Theorem:1.3} concerns with the existence of solutions to problem~\eqref{eq:P} in the critical case $p=p_*$. 
\begin{theorem}
\label{Theorem:1.3}
	Let $\theta\in(0,2]$ with $N>\theta$, $p=p_*$, and $T_*\in(0,\infty)$.
	Then there exists $\varepsilon>0$ such that 
	if $\mu\in\mathfrak L^{1,\infty}_{\rm ul}(\log\mathfrak L)^{(N-\theta)/\theta}$ satisfies
	\begin{equation}
	\label{eq:1.10}
		|||\mu|||_{1,(N-\theta)/\theta;T^{1/\theta}}\le\varepsilon\quad\mbox{for some $T\in(0,T_*]$},
	\end{equation}
	then problem~\eqref{eq:P} possesses a solution~$u$ in ${\mathbb R}^N\times[0,T)$ such that
	\begin{equation}
	\label{eq:1.11}
		\sup_{t\in(0,T)}\,\|u(t)\|_{p,N/\theta;T^{1/\theta}}\le C|||\mu|||_{1,(N-\theta)/\theta;T^{1/\theta}},
	\end{equation}
	where $C$ is a positive constant depending only on $T_*$, $N$, and $\theta$.
\end{theorem}
We remark that the constant $\varepsilon$ in Theorem~\ref{Theorem:1.3} tends to $0$ as $T_*\to\infty$ 
(see also Remark~\ref{Remark:1.1}).  
As a corollary of Theorem~\ref{Theorem:1.3}, we have:
\begin{corollary}
\label{Corollary:1.1}
	Let $\theta\in(0,2]$ with $N>\theta$ and $p=p_*$.
	Let $\mu\in {\mathfrak L}_{{\rm ul}}^{1,\infty}(\log{\mathfrak L})^\alpha$ for some $\alpha>(N-\theta)/\theta$.
	Then problem~\eqref{eq:P} possesses a local-in-time solution.
\end{corollary}
Furthermore, in the critical case $p=p_*$, it follows from \eqref{eq:1.2} and $N>\theta$ that 
\begin{equation}
\label{eq:1.12}
\mu_c^*(s)=O\left(s^{-1}|\log s|^{-\frac{N}{\theta}}\right),
\quad
\int_0^s \mu_c^*(\tau)\,d\tau=O\left(|\log s|^{-\frac{N-\theta}{\theta}}\right),\quad\mbox{as $s\to +0$}.
\end{equation}
These imply that $\mu_c\chi_{B(0,1)}\in {\mathfrak L}^{1,\infty}(\log{\mathfrak L})^{(N-\theta)/\theta}$. 
Then Theorem~\ref{Theorem:1.3} implies Theorem~\ref{Theorem:1.2}~(5), 
that is, we have:
\begin{corollary}
\label{Corollary:1.2}
	Let $\theta\in(0,2]$ with $N>\theta$ and $p=p_*$.
	Then there exists $\varepsilon>0$ such that if $\mu\in{\mathcal L}$ satisfies 
	$$
	|\mu(x)|\le\varepsilon\mu_c(x)+C, \quad\mbox{a.a.~$x\in{\mathbb R}^N$},
	$$
	for some $C\ge 0$, 
	then problem~\eqref{eq:P} possesses a local-in-time solution.
\end{corollary}

Next, we state our main theorem in the supercritical case $p>p_*$. 
\begin{theorem}
\label{Theorem:1.4}
	Let $\theta\in(0,2]$ with $N>\theta$ and $p>p_*$. 
	Then there exists $\varepsilon>0$ such that 
	if $\mu \in L^{r_*,\infty}_{\rm ul}$ satisfies
	$$
	\|\mu\|_{r_*,0;T^{1/\theta}} \le \varepsilon\,\,\,\,\mbox{for some $T\in(0,\infty]$, where $\displaystyle{r_*:=\frac{N(p-1)}{\theta p}}$},
	$$
	then the problem $(P)$ possesses a solution~$u$ in ${\mathbb R}^N\times[0,T)$ such that
	\begin{equation}
	\label{eq:1.13}
	\sup_{t\in(0,T)}\,\|u(t)\|_{pr_*,0;T^{1/\theta}} \le C \|\mu\|_{r_*,0;T^{1/\theta}},
	\end{equation}
	where $C$ is a positive constant depending only on $p, N$, and $\theta$.
\end{theorem}
Similarly to Corollary~\ref{Corollary:1.1}, we have:
\begin{corollary}
\label{Corollary:1.3}
	Let $\theta\in(0,2]$ with $N>\theta$ and $p>p_*$. 
	Let $\mu\in L_{{\rm ul}}^{r_*,\infty}(\log L)^\alpha$ for some $\alpha>0$. 
	Then problem~\eqref{eq:P} possesses a local-in-time solution.
\end{corollary}
In the supercritical case $p>p_*$, it follows from \eqref{eq:1.2} that 
\begin{equation}
\label{eq:1.14}
\mu_c^*(s)=O\left(s^{-\frac{\theta p}{N(p-1)}}\right)=O\left(s^{-\frac{1}{r_*}}\right)\quad\mbox{as $s\to +0$}.
\end{equation}
This implies that $\mu_c\chi_{B(0,1)}\in L^{r_*,\infty}$. 
Then Theorem~\ref{Theorem:1.4} implies Theorem~\ref{Theorem:1.2}~(6), 
that is, we have: 
\begin{corollary}
\label{Corollary:1.4}
	Let $\theta\in(0,2]$ with $N>\theta$ and $p>p_*$. 
	Then there exists $\varepsilon>0$ such that if $\mu\in{\mathcal L}$ satisfies
	$$
	|\mu(x)|\le\varepsilon \mu_c(x)+C, \quad\mbox{a.a.~$x\in{\mathbb R}^N$},
	$$
	for some $C>0$ {\rm({\it resp.~$C=0$})}, 
	then problem~\eqref{eq:P} possesses a local-in-time solution {\rm ({\it resp.~global-in-time solution})}.
\end{corollary}

Let $G$ be the fundamental solution to the fractional heat equation
\begin{equation}
\label{eq:1.15}
\partial_t v+(-\Delta)^{\frac{\theta}{2}}v=0\quad\mbox{in}\quad{\mathbb R}^N\times(0,\infty),
\end{equation}
where $\theta\in(0,2]$. 
For any $\varphi$ in $L^1_{{\rm loc}}$,
we set
\begin{equation}
\label{eq:1.16}
(S(t)\varphi)(x):=\int_{{\mathbb R}^N}G(x-y,t)\varphi(y)\,dy,
\quad (x,t)\in{\mathbb R}^N\times(0,\infty).
\end{equation}
The proofs of Theorems~\ref{Theorem:1.3} and \ref{Theorem:1.4} are based
 on decay estimates of  the fractional heat semigroup $S(t)$ in uniformly local Zygumnd spaces. 
Let $1\le r\le q\le\infty$ and $\alpha$, $\beta\in[0,\infty)$ with $\alpha\le\beta$ if $r=q$. 
In  \cite{IIK} Ioku {\it et al.} studied decay estimates of the fractional heat semigroup $S(t)$ 
in the spaces ${\mathfrak L}^{q,\infty}(\log{\mathfrak L})^\alpha$:
$$
S(t):{\mathfrak L}^{r,\infty}_{{\rm ul}}(\log{\mathfrak L})^\alpha\rightarrow {\mathfrak L}^{q,\infty}_{{\rm ul}}(\log{\mathfrak L})^{\beta},
$$
and established the following estimate
\begin{equation}
\label{eq:1.17}
|||S(t)\varphi|||_{q,\beta;T^{1/\theta}}\le Ct^{-\frac{N}{\theta}\left(\frac{1}{r}-\frac{1}{q}\right)}
\Phi(t^{-1})^{-\frac{\alpha}{r}+\frac{\beta}{q}}|||\varphi|||_{r,\alpha;T^{1/\theta}},\quad t\in(0,T),
\end{equation}
for $\varphi\in {\mathfrak L}^{r,\infty}_{{\rm ul}}(\log{\mathfrak L})^\alpha$ and $T\in(0,\infty]$
(see \cite{IIK}*{Proposition~3.2}).
This estimate played a crucial role to obtain sharp sufficient conditions for the existence of solutions 
to the critical semilinear fractional heat equation 
$$
\partial_t u+(-\Delta)^{\frac{\theta}{2}}u=u^{1+\frac{\theta}{N}},\quad x\in{\mathbb R}^N,\,\,\,t>0,
$$
with singular initial data. 
Recently, in \cite{FIK} Fujishima {\it et al.} studied decay estimates of the heat semigroup $e^{t\Delta}$, 
that is, $S(t)$ with $\theta=2$:
\begin{equation}
\label{eq:1.18}
\begin{array}{ll}
S(t): L^{r,\infty}_{{\rm ul}}(\log L)^\alpha\rightarrow{\mathfrak L}^{q,\infty}_{{\rm ul}}(\log{\mathfrak L})^{\beta}\subset L^{q,\infty}(\log L)^{\beta}
 & \quad\mbox{if}\quad 1<r<q\le\infty,\vspace{5pt}\\
S(t): L^{r,\infty}_{{\rm ul}}(\log L)^\alpha\rightarrow L^{q,\infty}_{{\rm ul}}(\log L)^{\beta}
 & \quad\mbox{if}\quad 1<r\le q\le\infty,
\end{array}
\end{equation}
and obtained the estimate
\begin{align}
&
\label{eq:1.19}
	|||S(t)\varphi|||_{q,\beta;T^{1/\theta}}\le Ct^{-\frac{N}{2}\left(\frac{1}{r}-\frac{1}{q}\right)}
	\Phi(t^{-1})^{-\frac{\alpha}{r}+\frac{\beta}{q}}\| \varphi \|_{r,\alpha;T^{1/\theta}},
	\quad 1<r<q \le \infty, \quad t\in(0,T),
	\\
&
\notag
	\|S(t)\varphi\|_{q,\beta;T^{1/\theta}}\le Ct^{-\frac{N}{2}\left(\frac{1}{r}-\frac{1}{q}\right)}
	\Phi(t^{-1})^{-\frac{\alpha}{r}+\frac{\beta}{q}}\| \varphi \|_{r,\alpha;T^{1/\theta}},
	\quad\,\,\,\, 1<r\le q \le \infty, \quad t\in(0,T)
\end{align}
(see \cite{FIK}*{Proposition~3.2}).
These estimates and their generalizations were useful in \cite{FIK} for the  study of sharp sufficient conditions 
for the existence of solutions to semilinear parabolic systems. 

In this paper, for the proof of Theorem \ref{Theorem:1.3}, 
we firstly establish the critical linear inhomogeneous estimate
\begin{equation}
	\label{eq:1.20}
	\sup_{0<t<T}\bigg\|\int_0^t S(t-s)\mu\,ds\,\bigg\|_{N/(N-\theta),N/\theta;T^{1/\theta}}
	\le C|||\mu|||_{1,(N-\theta)/\theta;T^{1/\theta}}
\end{equation}
for all $T \in (0,\infty]$. Note that this homogeneous estimate \eqref{eq:1.20} does not follow from
a direct calculation based on estimate \eqref{eq:1.17} since
 \[
	\| S(t) \mu \|_{N/(N-\theta), N/\theta;T^{1/\theta}}
	\le
	||| S(t) \mu |||_{N/(N-\theta), N/\theta;T^{1/\theta}}
	\le
	C t^{-1} |||\mu|||_{1,(N-\theta)/\theta;T^{1/\theta}},
	\quad t\in(0,T).
\]
In the proof, we employ the real interpolation method
developed through the analysis of the incompressible Navier-Stokes equations
with time-dependent external forces in the scaling critical framework
(see e.g., \cites{Yamazaki, Meyer, OT}). 
We prove a real interpolation relation
\begin{equation}
	\label{eq:1.21}
	({\mathfrak L}^{q_0,\infty}(\log{\mathfrak L})^{\alpha \frac{q_0}{q}},
	\, {\mathfrak L}^{q_1,\infty}(\log{\mathfrak L})^{\alpha \frac{q_1}{q}})_{\kappa, \infty}
	\hookrightarrow L^{q,\infty}(\log L)^\alpha
\end{equation}
for $1\le q_0<q<q_1<\infty$ satisfying
$1/q=(1-\kappa)/q_0+\kappa/q_1$
with $\kappa\in(0,1)$ and $\alpha\in[0,\infty)$.
This together with estimate~\eqref{eq:1.17}
enables us to show the critical linear inhomogeneous estimate~\eqref{eq:1.20}.
Next, we establish the critical nonlinear inhomogeneous estimates.
We shall show that
for any $T\in(0,\infty)$, there exists $C_T>0$ depending on $T$ such that
\begin{equation*}
	\sup_{0<t<T}\bigg\|\int_0^t S(t-s) |u(s)|^{p-1}u(s)  \,ds\,
	\bigg\|_{N/(N-\theta),N/\theta;T^{1/\theta}}
	\le C_T
	\sup_{0<t<T} \| u(t) \|_{N/(N-\theta),N/\theta;T^{1/\theta}}^p,
\end{equation*}
where $p=p_*=N/(N-\theta)$.
To this end, we treat the critical case $r=1$ in estimate~\eqref{eq:1.19},
which requires the restriction $r>1$,
and establish the following local estimates with modifications of the logarithmic order:
for $\alpha>0$ and $T\in(0,\infty)$, there exists $C_T>0$ depending on $T$ such that
\begin{equation}
	\label{eq:1.22}
	|||S(t)\varphi|||_{q,\beta;T^{1/\theta}}
	\le C_Tt^{-\frac{N}{\theta}\left(1-\frac{1}{q}\right)}\Phi(t^{-1})^{-\alpha+\frac{\beta}{q}}\|\varphi\|_{1,\alpha+1;T^{1/\theta}},\quad t\in(0,T)
\end{equation}
(see Proposition~\ref{Proposition:3.2}~(3)).
Using estimate \eqref{eq:1.22}, 
we apply the real interpolation relation \eqref{eq:1.21} again to prove that 
\[
	u\mapsto 
	\int_0^t S(t-s)|u(s)|^{p-1}u(s)\,ds
\]
is a contraction map in a neighborhood of the origin in $C((0,T):L^{p,\infty}_{{\rm ul}}(\log L)^{p\gamma})$,
where $\gamma:=(N-\theta)/\theta$. 
Then, by the contraction mapping theorem we find a solution~$u$ to problem~\eqref{eq:P} 
under assumption~\eqref{eq:1.10}.
Furthermore, we modify the arguments in the proof of Theorem~\ref{Theorem:1.3} to prove Theorem~\ref{Theorem:1.4}.
We remark that, in estimate \eqref{eq:1.22}, we cannot take $T=\infty$.
In fact, if estimate \eqref{eq:1.22} holds for $T=\infty$, then we have
$$
\|S(t)\varphi\|_{L^q(\log L)^\beta}
\le Ct^{-\frac{N}{\theta}\left(1-\frac{1}{q}\right)}\Phi(t^{-1})^{-\alpha+\frac{\beta}{q}}\|\varphi\|_{L^1(\log L)^{\alpha+1}},\quad t\in(0,\infty),
$$
and applying the same arguments as in the proof of Theorem~\ref{Theorem:1.3},
we can prove the existence of nontrivial global-in-time solutions to problem~\eqref{eq:P} with $p=p_*$
under assumption~\eqref{eq:1.10}.
This contradicts \eqref{eq:1.1}.
Therefore the condition $T<\infty$ is essential in estimate \eqref{eq:1.22}.
We also give an improvement of estimates \eqref{eq:1.17} and \eqref{eq:1.18} to
\begin{equation*}
		S(t): L^{r,\infty}_{{\rm ul}}(\log L)^\alpha\rightarrow
		L^{q}_{{\rm ul}}(\log L)^{\beta}
		\quad\mbox{if}\quad 1<r<q\le\infty
\end{equation*}
with 
	\begin{equation*}
		|S(t)\varphi|_{q,\beta;T^{1/\theta}}\le Ct^{-\frac{N}{2}\left(\frac{1}{r}-\frac{1}{q}\right)}
		\Phi(t^{-1})^{-\frac{\alpha}{r}+\frac{\beta}{q}}\| \varphi \|_{r,\alpha;T^{1/\theta}},
		\quad t\in(0,T).
	\end{equation*}
See Propositions \ref{Proposition:3.1} and \ref{Proposition:3.2}.
\medskip

The rest of this paper is organized as follows. 
In Section~2 we obtain preliminary lemmas on our Zygmund spaces (see Section~2.1). 
Furthermore, we recall Hardy's inequalities (see Section 2.2) and a decay estimate of 
the fractional heat semigroup $S(t)$ in Lebesgue spaces (see Section 2.3). 
In Section~3 we obtain decay estimates of $S(t)$ in our Zygmund spaces. 
In Sections~4 we employ the real interpolation method together with estimates of $S(t)$ given in Section~3. 
Furthermore, we apply the contraction mapping theorem to prove Theorem~\ref{Theorem:1.3}. 
We also prove Corollaries~\ref{Corollary:1.1} and~\ref{Corollary:1.2}.
In Section~5 we modify the arguments in Section~4, and prove Theorem~\ref{Theorem:1.4}, 
Corollary~\ref{Corollary:1.3}, and Corollary~\ref{Corollary:1.4}.
\section{Preliminaries}
In this section we recall some lemmas on Zygmund spaces and Hardy's inequalities. 
In all that follows we will use $C$ to denote generic positive constants and point out that $C$  
may take different values  within a calculation.
We often use the following properties of the function~$\Phi$ (see \eqref{eq:1.3}):
\begin{itemize}
  \item[($\Phi$1)] 
  the function: $(0,\infty)\ni s\mapsto\Phi(s^{-1})\in(1,\infty)$ is non-increasing;
  \item[($\Phi$2)] 
  for any fixed $k>0$, 
  $$
  \Phi(s^{-1})\le\ C\Phi(ks^{-1})\le C\Phi(s^{-k})\le C\Phi(s^{-1}),\quad s\in(0,\infty)
  $$
  (see e.g., \cite{FIK}*{Lemma~3.5~(1)});
  \item[($\Phi$3)]
  for any fixed $\alpha\in{\mathbb R}$ and $\delta>0$, 
  $$
  s_1^\delta \Phi(s_1^{-1})^\alpha\le Cs_2^\delta \Phi(s_2^{-1})^\alpha,
  \qquad
  s_1^{-\delta} \Phi(s_1^{-1})^\alpha\ge Cs_2^{-\delta} \Phi(s_2^{-1})^\alpha,
  $$
  for $s_1$, $s_2\in(0,\infty)$ with $s_1\le s_2$ (see e.g., \cite{FIK}*{Lemma~3.5~(2)}).
\end{itemize} 
\subsection{Preliminaries for Zygmund spaces}
We first recall some properties of non-increasing rearrangements of $f$, $g\in{\mathcal L}$ 
(see e.g., \cite{Grafakos}*{Proposition~1.4.5}).
\begin{itemize}
\item[(R1)]
For any  $q\in(0,\infty)$ and $k\in{\mathbb R}$,
$$
	(kf)^*=|k|f^*,\quad (|f|^q)^*=(f^*)^q,\quad 
	\int_{{\mathbb R}^N}|f(x)|^q\,dx=\int_0^\infty f^*(s)^q\,ds,\quad
	f^*(0)=\|f\|_{L^\infty}.
$$
\item[(R2)]
$(f+g)^*(t+s)\le f^*(t)+g^*(s)$ and $(fg)^*(t+s)\le f^*(t) g^*(s)$ for $t$, $s\in[0,\infty)$.
\item[(R3)]
For any $E\in{\mathcal M}$ with $|E|\not=0$, 
$(\chi_E)^*(s)=\chi_{[0,|E|)}(s)$ for $s\in[0,\infty)$.
\item[(R4)]
$\displaystyle{f^{**}(s):=\frac{1}{s}\int_0^s f^*(\tau)\,d\tau\ge f^*(s)}$ for $s\in(0,\infty)$.
\end{itemize}
In addition, we have:
\begin{itemize}
  \item[(R5)] 
  $\displaystyle{f^{**}(s)=\frac{1}{s}\sup_{|E|=s}\int_E |f(x)|\,dx}$ for $s\in(0,\infty)$ 
  (see \cite{BS}*{Chapter 2, Proposition 3.3});
  \item[(R6)]
  {(O'Neil's inequality)} 
  $$
  (f*g)^{**}(s)\le\int_s^\infty  f^{**}(\tau)g^{**}(\tau)\,d\tau
  $$ 
  for $s\in(0,\infty)$, where 
  $\displaystyle{(f*g)(x)=\int_{{\mathbb R}^N}f(x-y)g(y)\,dy}$ 
  (see \cite{ONeil}*{Lemma~1.6});
  \item[(R7)]
  $\displaystyle{(fg)^{**}(s)\le \frac{1}{s}\int_0^s f^*(\tau)g^*(\tau)\,d\tau}$ for $s\in(0,\infty)$ 
  (see \cite{ONeil}*{Theorem 3.3}).
\end{itemize}

Next, we recall some properties of Zygmund spaces. 
Let $q\in[1,\infty)$ and $\alpha\ge 0$.
For any $f\in{\mathcal L}$, 
by  (R1), (R5), and \eqref{eq:1.6}
we have
\begin{equation}
\label{eq:2.1}
\begin{split}
	\|f\|_{{\mathfrak L}^{q,\infty}(\log{\mathfrak L})^\alpha}
	& 
	=\sup_{s>0}\,\left\{\Phi(s^{-1})^\alpha\int_0^s (|f|^q)^*(\tau)\,d\tau\right\}^{\frac{1}{q}}
	\\
	& 
	=\sup_{s>0}\,\left\{s\Phi(s^{-1})^\alpha (|f|^q)^{**}(s)\right\}^{\frac{1}{q}}
	\\
 	& 
	=\sup_{s>0}\,\left\{\Phi(s^{-1})^\alpha\sup_{|E|=s}\int_E |f(x)|^q\,dx\right\}^{\frac{1}{q}}.
\end{split}
\end{equation}
On the other hand, 
for any $E\in{\mathcal M}$ with $|E|<\infty$, 
it follows from (R3) and (R7) that 
\begin{equation}
\begin{split}
\label{eq:2.2}
	(f\chi_E)^{**}(s) & =(f\chi_E\chi_E)^{**}(s)\\
	 & \le\frac{1}{s}\int_0^s (f\chi_E)^*(\tau)(\chi_E)^*(\tau)\,d\tau
	 =\frac{1}{s}\int_0^{\min\{s,|E|\}} (f\chi_E)^*(\tau)\,d\tau,\quad s\in(0,\infty).
\end{split}
\end{equation}
For any $\beta\in[\alpha,\infty)$, 
by ($\Phi1$), \eqref{eq:2.1}, and \eqref{eq:2.2} we have
\begin{align*}
	\|f\chi_E\|_{{\mathfrak L}^{1,\infty}(\log{\mathfrak L})^\alpha}
	& 
	=\sup_{s>0}\,\left\{s\Phi(s^{-1})^\alpha  (f\chi_E)^{**}(s)\right\}\\
	&
	\le \sup_{s>0}\,\left\{\Phi(s^{-1})^\alpha \int_0^{\min\{s,|E|\}} (f\chi_E)^*(\tau)\,d\tau\right\}\\
 	& 
 	=\sup_{0<s\le |E|}\left\{\Phi(s^{-1})^{\beta+\alpha-\beta}
	\int_0^{\min\{s,|E|\}} (f\chi_E)^*(\tau)\,d\tau\right\}\\
 	& 
 	\le\Phi(|E|^{-1})^{\alpha-\beta}
 	\sup_{0<s\le |E|}\left\{\Phi(s^{-1})^{\beta}
	\int_0^{\min\{s,|E|\}} (f\chi_E)^*(\tau)\,d\tau\right\}\\
	& 
 	\le\Phi(|E|^{-1})^{\alpha-\beta}
 	\sup_{s>0}\left\{\Phi(s^{-1})^{\beta}
	\int_0^s (f\chi_E)^*(\tau)\,d\tau\right\}\\
 	& 
 	=\Phi(|E|^{-1})^{\alpha-\beta}
	\|f\chi_E\|_{{\mathfrak L}^{1,\infty}(\log{\mathfrak L})^\beta}.
\end{align*}
This together with ($\Phi$2) implies that
\begin{equation}
\label{eq:2.3}
	|||f|||_{1,\alpha;\rho}\le C\Phi(\rho^{-1})^{\alpha-\beta}|||f|||_{1,\beta;\rho}
\end{equation}
for $0\le\alpha\le\beta$ and $\rho\in(0,\infty]$. 
%
\begin{lemma}
\label{Lemma:2.1}
	Let $q$, $q_1$, $q_2\in[1,\infty)$ and $\alpha$, $\alpha_1$, $\alpha_2\in[0,\infty)$ be such that 
	\begin{equation}
	\label{eq:2.4}
		\frac{1}{q}=\frac{1}{q_1}+\frac{1}{q_2},\quad 
		\frac{\alpha}{q}=\frac{\alpha_1}{q_1}+\frac{\alpha_2}{q_2}.
	\end{equation}
	Then there exists $C>0$, independent of $q_1$, $q_2$, $\alpha$, $\alpha_1$, and $\alpha_2$, 
	such that
	\begin{align}
	\label{eq:2.5}
	 \|fg\|_{L^{q,\infty}(\log L)^\alpha}
	 & \le C \|f\|_{L^{q_1,\infty}(\log L)^{\alpha_1}}\|g\|_{L^{q_2,\infty}(\log L)^{\alpha_2}},\\
	\label{eq:2.6}
	\|fg\|_{q,\alpha;\rho}
	 & \le C \|f\|_{q_1,\alpha_1;\rho} \|g\|_{q_2,\alpha_2;\rho},
	\end{align}
	for $f$, $g\in{\mathcal L}$ and $\rho\in(0,\infty]$. 
\end{lemma}
{\bf Proof.}
Let $q$, $q_1$, $q_2\in[1,\infty)$ and $\alpha$, $\alpha_1$, $\alpha_2\ge 0$, and assume \eqref{eq:2.4}. 
Since 
$$
	\Phi(s^{-1})\le
	\Phi\left(\left(\frac{s}{2}\right)^{-1}\right), \quad s>0,
$$
it follows from (R2), and \eqref{eq:1.5} that
\begin{align*}
	\|fg\|_{L^{q,\infty}(\log L)^\alpha}
 	& 
	=\sup_{s>0}\,\left\{s^\frac{1}{q}\Phi(s^{-1})^\frac{\alpha}{q}  (fg)^{*}(s)\right\}
	\\
 	& 
	\le\sup_{s>0}\,\left\{s^\frac{1}{q} \Phi(s^{-1})^\frac{\alpha}{q} f^*\left(\frac{s}{2}\right)g^*\left(\frac{s}{2}\right)\right\}
	\\
 	& 
	\le 
	 2^{\frac{1}{q}}
	 \sup_{s>0}\,\biggr\{\Phi\left(\left(\frac{s}{2}\right)^{-1}\right)^{\frac{\alpha_1}{q_1}+\frac{\alpha_2}{q_2}}
	\left(\frac{s}{2}\right)^{\frac{1}{q_1}+\frac{1}{q_2}}f^*\left(\frac{s}{2}\right)g^*\left(\frac{s}{2}\right)\biggr\}
	\\
 	& 
	\le C
	\sup_{s>0}\,\left\{s\Phi(s^{-1})^{\alpha_1}f^*(s)^{q_1}\right\}^\frac{1}{q_1}\cdot
 	\sup_{s>0}\,\left\{s\Phi(s^{-1})^{\alpha_2}g^*(s)^{q_2}\right\}^\frac{1}{q_2}
	\\
 	& 
	=C
	\|f\|_{L^{q_1,\infty}(\log L)^{\alpha_1}}\|g\|_{L^{q_2,\infty}(\log L)^{\alpha_2}},
	\quad f,g\in{\mathcal L}.
\end{align*}
This implies \eqref{eq:2.5}. 
Furthermore,
by \eqref{eq:2.5} we have 
\begin{align*}
	\|fg\|_{q,\alpha;\rho}
 	& 
	=\sup_{x\in{\mathbb R}^N}\|fg\chi_{B(x,\rho)}\|_{L^{q,\infty}(\log L)^\alpha}
	\\
 	& 
	\le C
	\sup_{x\in{\mathbb R}^N}\left\{\|f\chi_{B(x,\rho)}\|_{L^{q_1,\infty}(\log L)^{\alpha_1}}
	\|g\chi_{B(x,\rho)}\|_{L^{q_2,\infty}(\log L)^{\alpha_2}}\right\}
	\\
	 & 
	 \le C
	 \sup_{x\in{\mathbb R}^N}\|f\chi_{B(x,\rho)}\|_{L^{q_1,\infty}(\log L)^{\alpha_1}}\cdot 
  	\sup_{x\in{\mathbb R}^N}\|g\chi_{B(x,\rho)}\|_{L^{q_2,\infty}(\log L)^{\alpha_2}}
	\\
 	& 
	=C
	\|f\|_{q_1,\alpha_1;\rho}\|g\|_{q_2,\alpha_2;\rho},
	\quad f,g\in{\mathcal L},\,\,\,\rho\in(0,\infty]. 
\end{align*}
This implies \eqref{eq:2.6}. Thus Lemma~\ref{Lemma:2.1} follows. 
$\Box$
\begin{lemma}
\label{Lemma:2.2}
	Let $q\in[1,\infty)$, $\alpha\in[0,\infty)$, and $\rho\in(0,\infty]$. 
	Then 
	$$
	\||f|^r\|_{L^{q,\infty}(\log L)^\alpha}=\|f\|_{L^{rq,\infty}(\log L)^{\alpha}}^r,
	\quad 
	\||f|^r\|_{q,\alpha;\rho}=\|f\|_{rq,\alpha;\rho}^r,
	$$
	for $f\in{\mathcal L}$ and $r\in(0,\infty)$ with $rq\ge 1$. 
\end{lemma}
{\bf Proof.}
Let $f\in{\mathcal L}$, $\alpha\in[0,\infty)$, and $\rho\in(0,\infty]$.
It follows from (R1) and \eqref{eq:1.5} that 
\begin{align*}
	\||f|^r\|_{L^{q,\infty}(\log L)^\alpha} 
	& 
	=\sup_{s>0}\,\left\{s\Phi(s^{-1})^\alpha (|f|^r)^*(s)^q\right\}^{\frac{1}{q}}
	\\
	& 
	=\sup_{s>0}\,\left\{s\Phi(s^{-1})^{\alpha} f^*(s)^{rq}\right\}^{\frac{r}{rq}}
	 =\|f\|_{L^{rq,\infty}(\log L)^{\alpha}}^r.
\end{align*}
This implies that
$$
	\||f|^r\|_{q,\alpha;\rho}
	=\sup_{x\in{\mathbb R}^N}\| |f|^r\chi_{B(x,\rho)}^r\|_{L^{q,\infty}(\log L)^\alpha}
	=\sup_{x\in{\mathbb R}^N}\| |f|\chi_{B(x,\rho)}\|_{L^{rq,\infty}(\log L)^{\alpha}}^r
	=\|f\|_{rq,\alpha;\rho}^r
$$
for $\rho\in(0,\infty]$. 
Thus Lemma~\ref{Lemma:2.2} follows. 
$\Box$
\subsection{Hardy's inequalities}
We recall the following two lemmas on Hardy's inequality.
(See \cite{Muckenhoup}*{Theorems~1 and 2} and \cite{OK}*{Theorems~5.9 and 6.2}.)
Throughout this paper, for any $q\in[1,\infty]$, we denote by $q'$ the H\"older conjugate of $q$, that is, 
$q'=q/(q-1)$ if $q\in(1,\infty)$, $q'=\infty$ if $q=1$, and $q'=1$ if $q=\infty$.
\begin{lemma}
\label{Lemma:2.3}
Let $q\in[1,\infty]$ and $0\le a<b\le\infty$.
Let $U$ and $V$ be locally integrable functions in $[a,b)$. 
Then there exists $c_*>0$ such that 
$$
\|U\tilde{F}\|_{L^q((a,b))}
\le c_*\|Vf\|_{L^q((a,b))}
\quad\mbox{with}\quad
\tilde{F}(s):=\int_a^s f(\tau)\,d\tau
$$
holds for locally integrable functions $f$ in $[a,b)$ 
if and only if 
$$
B_1:=
\sup_{a<s<b}\,
\left\{\|U\|_{L^q((s,b))}\|V^{-1}\|_{L^{q'}((a,s))}\right\}<\infty.
$$
Here $c_*\le q^{1/q}(q')^{1/q'}B_1$ for $q\in(1,\infty)$ and $c_*=B_1$ if $q=1,\infty$.
\end{lemma}
\begin{lemma}
\label{Lemma:2.4}
Let $q\in[1,\infty]$ and $0\le a<b\le\infty$.
Let $U$ and $V$ be locally integrable functions in $[a,b)$. 
Then there exists $c_{**}>0$ such that 
$$
\|U\hat{F}\|_{L^q((a,b))}
\le c_{**}\|Vf\|_{L^q((a,b))}
\quad\mbox{with}\quad
\hat{F}(s):=\int_s^b f(\tau)\,d\tau
$$
holds for locally integrable functions $f$ in $(a,b)$ with $f\in L^1((c,b))$ for any $c\in(a,b)$ 
if and only if 
$$
B_2:=\sup_{a<s<b}\left\{\|U\|_{L^q((a,s))}\|V^{-1}\|_{L^{q'}((s,b))}\right\}<\infty.
$$
Here $c_{**}\le q^{1/q}(q')^{1/q'}B_2$ for $q\in(1,\infty)$ and $c_{**}=B_2$ if $q=1,\infty$.
\end{lemma}
\subsection{Fundamental solutions}
Let $\theta\in(0,2]$. 
Let $G$ be the fundamental solution to fractional heat equation~\eqref{eq:1.15}. 
The function $G$ is positive and smooth in ${\mathbb R}^N\times(0,\infty)$ 
and it satisfies
\begin{equation}
\label{eq:2.7}
\begin{array}{ll}
\displaystyle{G(x,t)=(4\pi t)^{-\frac{N}{2}}\exp\left(-\frac{|x|^2}{4t}\right)\le Ch_t(x)} & \quad\mbox{if}\quad \theta=2,\vspace{5pt}\\
C^{-1}h_t(x)\le G(x,t)\le Ch_t(x) & \quad\mbox{if}\quad 0<\theta<2,
\end{array}
\end{equation}
for $(x,t)\in{\mathbb R}^N\times(0,\infty)$ (see e.g., \cites{BJ, BraK, S}). 
Here
$h_t(x):=t^{-\frac{N}{\theta}}\left(1+t^{-\frac{1}{\theta}}|x|\right)^{-N-\theta}$. 
Furthermore, $G$ satisfies 
\begin{eqnarray*}
 & \bullet & G(x,t)=t^{-\frac{N}{\theta}}G\left(t^{-\frac{1}{\theta}}x,1\right),
 \quad \int_{{\mathbb R}^N}G(x,t)\,dx=1,\\
 & \bullet &  \mbox{$G(\cdot,1)$ is radially symmetric and $G(x,1)\le G(y,1)$ if $|x|\ge |y|$},\\
 & \bullet & G(x,t)=\int_{{\mathbb R}^N}G(x-y,t-s)G(y,s)\,dy,
\end{eqnarray*}
for $x$, $y\in{\mathbb R}^N$ and $0<s<t$.
Then it follows from Young's inequality that 
\begin{equation}
\label{eq:2.8}
\|S(t)\varphi\|_{L^q}\le Ct^{-\frac{N}{\theta}\left(\frac{1}{r}-\frac{1}{q}\right)}\|\varphi\|_{L^r}
\end{equation}
for $\varphi\in L^r$, $1\le r\le q\le\infty$, and $t>0$. 
(See \eqref{eq:1.16} for the definition of $S(t)\varphi$.)
\section{Decay estimates of $S(t)\varphi$}
In this section we prove the following proposition on decay estimates of $S(t)\varphi$ in Zygmund spaces. 
\begin{proposition}
\label{Proposition:3.1}
Let $\theta\in(0,2]$, $r$, $q\in[1,\infty]$ with $r<q$, and $\alpha$, $\beta\in[0,\infty)$. 
\begin{itemize}
\item[\rm (i)]
There exists $C_1>0$ such that
\begin{equation*}
\begin{split}
 & \|S(t)\varphi\|_{L^q(\log L)^\beta}+\|S(t)\varphi\|_{L^{q,\infty}(\log L)^\beta}+\|S(t)\varphi\|_{{\mathfrak L}^{q,\infty}(\log{\mathfrak L})^\beta}\\
 & \le C_1t^{-\frac{N}{\theta}\left(\frac{1}{r}-\frac{1}{q}\right)}\Phi(t^{-1})^{-\frac{\alpha}{r}+\frac{\beta}{q}}\|\varphi\|_{{\mathfrak L}^{r,\infty}(\log{\mathfrak L})^\alpha}
\end{split}
\end{equation*}
for $\varphi\in {\mathfrak L}^{r,\infty}(\log{\mathfrak L})^\alpha$ and $t>0$.
\item[\rm (ii)]
If $r>1$, then there exists $C_2>0$ such that
\begin{equation*}
\begin{split}
 & \|S(t)\varphi\|_{L^q(\log L)^\beta}+\|S(t)\varphi\|_{L^{q,\infty}(\log L)^\beta}+\|S(t)\varphi\|_{{\mathfrak L}^{q,\infty}(\log{\mathfrak L})^\beta}\\
 & \le C_2t^{-\frac{N}{\theta}\left(\frac{1}{r}-\frac{1}{q}\right)}\Phi(t^{-1})^{-\frac{\alpha}{r}+\frac{\beta}{q}}\|\varphi\|_{L^{r,\infty}(\log L)^\alpha}
\end{split}
\end{equation*}
for $\varphi\in L^{r,\infty}(\log L)^\alpha$ and $t>0$.
\item[\rm (iii)]
Let $E\in\mathcal M$ with $|E|<\infty$.
If $\alpha>0$, then there exists $C_3>0$ such that
\begin{equation*}
\begin{split}
 & \|S(t)(\varphi\chi_E)\|_{L^q(\log L)^\beta}+\|S(t)(\varphi\chi_E)\|_{L^{q,\infty}(\log L)^\beta}+\|S(t)(\varphi\chi_E)\|_{{\mathfrak L}^{q,\infty}(\log{\mathfrak L})^\beta}\\
 & \le C_3t^{-\frac{N}{\theta}\left(1-\frac{1}{q}\right)}\Phi(t^{-1})^{-\alpha+\frac{\beta}{q}}\|\varphi\chi_E\|_{L^{1,\infty}(\log L)^{\alpha+1}}
\end{split}
\end{equation*}
for $\varphi\in L^{1,\infty}(\log L)^{\alpha+1}$ and $t>0$. 
\end{itemize}
\end{proposition}
For the proof, 
we recall the following two lemmas (see \cite{IIK}*{Lemmas~3.1 and ~3.2}).
\begin{lemma}
\label{Lemma:3.1}
	{\rm (1)} 
  	For any $q>-1$ and $\alpha\in{\mathbb R}$, 
	there exists $C_1>0$ such that 
	$$
  	\int_0^s \tau^q\Phi(\tau^{-1})^\alpha\,d\tau
	\le C_1s^{q+1}\Phi(s^{-1})^\alpha,\quad s>0.
  	$$
	{\rm (2)}
	For any $\alpha<-1$ and $S\in(0,\infty)$, 
  	there exists $C_2>0$ such that 
  	$$
  	\int_0^s \tau^{-1}\Phi(\tau^{-1})^\alpha\,d\tau
	\le C_2\Phi(s^{-1})^{\alpha+1},
	\quad s\in(0,S).
  	$$
	\newline
	{\rm (3)}
  	For any $q<-1$ and $\alpha\in{\mathbb R}$, 
	there exists $C_3>0$ such that 
  	$$
  	\int_s^\infty \tau^q\Phi(\tau^{-1})^\alpha\,d\tau
	\le C_3s^{q+1}\Phi(s^{-1})^\alpha,
	\quad s>0.
  	$$
\end{lemma}
\begin{lemma}
\label{Lemma:3.2}
	Let $1\le r\le q<\infty$ and $\gamma\in{\mathbb R}$.
	Assume that $\gamma\ge 0$ if $r=q$. 
	Then there exists $C>0$ such that
	$$
		\int_0^\infty\tau^{q\left(1-\frac{1}{r}\right)}
		\Phi(\tau^{-1})^\gamma(h_t^*(\tau))^q\,d\tau
		\le Ct^{-\frac{Nq}{\theta}\left(\frac{1}{r}-\frac{1}{q}\right)}
		\Phi(t^{-1})^\gamma,\quad t>0.
	$$
	Here $h_t$ is as in Section~{\rm 2.3}. 
\end{lemma}
Next, we prove the following lemma. 
\begin{lemma}
\label{Lemma:3.3}
	{\rm (1)} 
  	For any $r\in[1,\infty)$ and $\alpha\in[0,\infty)$, 
	\begin{equation}
	\label{eq:3.1}
		\sup_{s>0}\,\bigg\{s\Phi(s^{-1})^{\alpha}
		f^{**}(s)^r\bigg\}^{\frac{1}{r}}
		\le \|f\|_{{\mathfrak L}^{r,\infty}(\log{\mathfrak L})^\alpha},
		\quad f\in{\mathcal L}.
	\end{equation}
	{\rm (2)}
  	For any $r\in(1,\infty)$ and $\alpha\in[0,\infty)$, there exists $C_1>0$ such that
	\begin{equation}
	\label{eq:3.2}
		\sup_{s>0}\,\bigg\{s\Phi(s^{-1})^{\alpha}
		f^{**}(s)^r\bigg\}^{\frac{1}{r}}
		\le C_1\|f\|_{L^{r,\infty}(\log L)^\alpha}, 
		\quad f\in{\mathcal L}.
	\end{equation}
	{\rm (3)}
  	For any $E\in\mathcal M$ with $|E|<\infty$ and $\alpha\in(0,\infty)$,
	there exists $C_2>0$ such that
	\begin{equation}
	\label{eq:3.3}
		\sup_{s>0}\,\bigg\{s\Phi(s^{-1})^{\alpha}(f\chi_E)^{**}(s)\bigg\}
		\le C_2\|f\chi_E\|_{L^{1,\infty}(\log L)^{\alpha+1}},
		\quad f\in{\mathcal L}.
	\end{equation}
	Furthermore, 
	$$
	\|f\chi_E\|_{{\mathfrak L}^{1,\infty}(\log{\mathfrak L})^\alpha}
	\le C_2\|f\chi_E\|_{L^{1,\infty}(\log L)^{\alpha+1}},\quad f\in{\mathcal L}.
	$$
\end{lemma}
{\bf Proof.}
Let $f\in{\mathcal L}$. 
We prove assertion (1). Let $r\in[1,\infty)$ and $\alpha\in[0,\infty)$. 
It follows from Jensen's inequality and \eqref{eq:1.6} that 
$$
	\sup_{s>0}\left\{s\Phi(s^{-1})^{\alpha} f^{**}(s)^r\right\}
	\le \sup_{s>0}\left\{s\Phi(s^{-1})^{\alpha}\cdot \frac{1}{s}\int_0^s f^*(\tau)^r\,d\tau\right\}
	=\|f\|_{{\mathfrak L}^{r,\infty}(\log{\mathfrak L})^\alpha}^r,
$$
which implies \eqref{eq:3.1}. Assertion~(1) follows. 

We prove assertion~(2) using Hardy's inequality. 
Let $r\in(1,\infty)$ and $\alpha\in[0,\infty)$. 
Set $U_1(\tau):=\tau^{\frac{1}{r}-1}\Phi(\tau^{-1})^{\frac{\alpha}{r}}$ and 
$V_1(\tau):=\tau^{\frac{1}{r}}\Phi(\tau^{-1})^{\frac{\alpha}{r}}$ for $\tau>0$.
Since $U_1$ is non-increasing by $r>1$, 
we have
$$
	\|U_	1\|_{L^\infty((s,\infty))}
	\le s^{\frac{1}{r}-1}\Phi(s^{-1})^{\frac{\alpha}{r}},\qquad s>0.
$$
This together with Lemma~\ref{Lemma:3.1}~(1) implies that
\begin{equation}
\label{eq:3.4}
	 \sup_{s>0}\,\left\{\|U_1\|_{L^\infty((s,\infty))}\int_0^s |V(\tau)|^{-1}\,d\tau\right\}
	\le \sup_{s>0}\left\{s^{\frac{1}{r}-1}\Phi(s^{-1})^{\frac{\alpha}{r}}
	\cdot Cs^{1-\frac{1}{r}}\Phi(s^{-1})^{-\frac{\alpha}{r}}\right\}<\infty.
\end{equation}
Applying Lemma~\ref{Lemma:2.3}, 
by \eqref{eq:1.5} and \eqref{eq:3.4} we obtain
\begin{align*}
\sup_{s>0}\, \left\{s\Phi(s^{-1})^{\alpha} f^{**}(s)^r\right\}^{\frac{1}{r}}
 & =\sup_{s>0}\, \left\{U_1(s)\int_0^sf^*(\tau)\,d\tau\right\}\\
 & \le C\sup_{s>0}\left\{V_1(s)f^*(s)\right\}
=C\|f\|_{L^{r,\infty}(\log L)^\alpha}.
\end{align*}
This implies \eqref{eq:3.2}, and assertion~(2) follows. 

We prove assertion~(3).
Let $E\in\mathcal M$ with $|E|<\infty$ and $\alpha\in(0,\infty)$. 
It follows from \eqref{eq:2.2} that
$$
\sup_{s>0}\,\left\{s\Phi(s^{-1})^\alpha (f\chi_E)^{**}(s)\right\}
\le \sup_{s>0}\,\left\{\Phi(s^{-1})^\alpha \int_0^{\min\{s,|E|\}} (f\chi_E)^*(\tau)\,d\tau\right\},
$$
which together with $\alpha\in(0,\infty)$ and ($\Phi1$) implies that
\begin{equation}
\begin{split}
\label{eq:3.5}
\sup_{s>0}\,\left\{s\Phi(s^{-1})^\alpha (f\chi_E)^{**}(s)\right\}
 & \le\sup_{0<s\le |E|}\left\{\Phi(s^{-1})^{\alpha}\int_0^{\min\{s,|E|\}} (f\chi_E)^*(\tau)\,d\tau\right\}\\
 & \le\sup_{0<s\le |E|}\left\{\Phi(s^{-1})^{\alpha}\int_0^s (f\chi_E)^*(\tau)\,d\tau\right\}.
\end{split}
\end{equation}
Set $U_2(\tau):=\Phi(\tau^{-1})^{\alpha}$ and $V_2(\tau):=\tau\Phi(\tau^{-1})^{\alpha+1}$ 
for $\tau>0$.
By Lemma~\ref{Lemma:3.1}~(2) we have
$$
	\int_0^s|V_2(\tau)|^{-1}\,d\tau\le C\Phi(s^{-1})^{-\alpha}\quad
	\mbox{for $s\in(0,|E|)$}.
$$
This together with $\alpha>0$ and ($\Phi1$) implies that
\begin{equation}
\label{eq:3.6}
\sup_{0<s<|E|}\bigg\{\|U_2\|_{L^\infty((s,|E|))}\int_0^s|V_2(\tau)|^{-1}\,d\tau\bigg\}
\le \sup_{0<s<|E|}\bigg\{\Phi(s^{-1})^\alpha\cdot C\Phi(s^{-1})^{-\alpha}\bigg\}<\infty.
\end{equation}
Applying Lemma~\ref{Lemma:2.3} to \eqref{eq:3.5},
by \eqref{eq:1.5} and \eqref{eq:3.6} we obtain 
\begin{align*}
	\sup_{s>0}\,\left\{s\Phi(s^{-1})^\alpha (f\chi_E)^{**}(s)\right\}
	& \le\sup_{0<s\le |E|}\left\{U_2(s)\int_0^s (f\chi_E)^*(\tau)\,d\tau\right\}\\
 	& \le C\sup_{0<s<|E|}\bigg\{V_2(s)(f\chi_E)^*(s)\bigg\}
	\le C\|f\chi_E\|_{L^{1,\infty}(\log L)^{\alpha+1}}. 
\end{align*}
This implies \eqref{eq:3.3}, and assertion~(3) follows. 
The proof of Lemma~\ref{Lemma:3.3} is complete.
$\Box$
\vspace{5pt}

Now we are ready to prove Proposition~\ref{Proposition:3.1}.
\vspace{5pt}
\newline
{\bf Proof of Proposition~\ref{Proposition:3.1}.}
The proof is divided into the following two cases: 
$$
1\le r<q<\infty,\qquad 1\le r< q=\infty.
$$ 
\underline{Step 1.}
Consider the case of $1\le r<q<\infty$. 
It follows from (R4), (R6), \eqref{eq:1.4}, and \eqref{eq:2.7} that 
\begin{equation}
\label{eq:3.7}
	\begin{aligned}
	\left\|S(t)\varphi\right\|_{L^q(\log L)^\beta}^q
 	& 
	=\int_0^\infty \Phi(\tau^{-1})^{\beta}
	\left((S(t)\varphi)^*(\tau)\right)^q\,d\tau
	\\
 	& 
	\le \int_0^\infty \Phi(\tau^{-1})^{\beta}
	(\left(S(t)\varphi\right)^{**}(\tau))^q\,d\tau
	\\
	 & 
	 \le C\int_0^\infty \Phi(\tau^{-1})^{\beta}
	\left(\int_\tau^\infty h_t^{**}(\eta)\varphi^{**}(\eta)\,d\eta\right)^q\,d\tau
	 \\
	 & 
	 = C\int_0^\infty \left(\Phi(\tau^{-1})^{\frac{\beta}{q}}
	 \int_\tau^\infty h_t^{**}(\eta)\varphi^{**}(\eta)\,d\eta\right)^q\,d\tau,
	 \quad t>0.
	\end{aligned}
\end{equation}
Set
$U_1(\tau):=\Phi(\tau^{-1})^{\frac{\beta}{q}}$ and $V_1(\tau):=\tau\Phi(\tau^{-1})^{\frac{\beta}{q}}$  
for $\tau>0$. 
It follows from Lemma~\ref{Lemma:3.1}~(1), (3) that
$$
	\sup_{s>0}\,\left(\int_0^s |U_1(\tau)|^q\,d\tau\right)^{\frac{1}{q}}
	\left(\int_s^\infty |V_1(\tau)|^{-q'}\,d\tau\right)^{\frac{1}{q'}}
	\le \sup_{s>0}\,\left\{Cs^{\frac{1}{q}}\Phi(s^{-1})^{\frac{\beta}{q}}\cdot 
	Cs^{-1+\frac{1}{q'}}\Phi(s^{-1})^{-{\frac{\beta}{q}}}\right\}<\infty.
$$
Then, by Lemma~\ref{Lemma:2.4} and \eqref{eq:3.7} we have 
\begin{equation}
\label{eq:3.8}
	\begin{aligned}
	& \left\|S(t)\varphi\right\|_{L^q(\log L)^\beta}^q\\
	& \le C\int_0^\infty \left(U_1(\tau)
	 \int_\tau^\infty h_t^{**}(\eta)\varphi^{**}(\eta)\,d\eta\right)^q\,d\tau\\
	 &  \le C\int_0^\infty \left(V_1(\tau)h_t^{**}(\tau)\varphi^{**}(\tau)\right)^q\,d\tau
 	=C\int_0^\infty \left(\tau\Phi(\tau^{-1})^{\frac{\beta}{q}} 
	h_t^{**}(\tau)\varphi^{**}(\tau)\right)^q\,d\tau
	\\
 	& 
	\le C\sup_{s>0}\,\left\{s\Phi(\tau^{-1})^{\alpha}
	(\varphi^{**}(s))^r\right\}^{\frac{q}{r}}\int_0^\infty\bigg( \tau^{1-\frac{1}{r}}
	\Phi(\tau^{-1})^{-\frac{\alpha}{r}+\frac{\beta}{q}}
	h_t^{**}(\tau)\bigg)^q\,d\tau
	\\
 	& 
	= C\sup_{s>0}\,\left\{s\Phi(s^{-1})^{\alpha}
	(\varphi^{**}(s))^r\right\}^{\frac{q}{r}}
	\int_0^\infty\bigg( \tau^{-\frac{1}{r}}
	\Phi(\tau^{-1})^\sigma
	\int_0^\tau h_t^*(\xi)\,d\xi\bigg)^q\,d\tau,\quad t>0,
	\end{aligned}
\end{equation}
where $\sigma:=-\alpha/r+\beta/q$. 
Set 
$U_2(\tau)=\tau^{-\frac{1}{r}}\Phi(\tau^{-1})^\sigma$ and $V_2(\tau)=\tau^{1-\frac{1}{r}}\Phi(\tau^{-1})^\sigma$ for $\tau>0$.
Since $q>r$ and $q'<r'$, 
by Lemma~\ref{Lemma:3.1}~(1), (3) we have 
\begin{equation}
\label{eq:3.9}
	\begin{aligned}
 	& 
	\sup_{s>0}\,\left(\int_s^\infty |U_2(\tau)|^q\,d\tau\right)^{\frac{1}{q}}
	\left(\int_0^s |V_2(\tau)|^{-q'}\,{d\tau}\right)^{\frac{1}{q'}}\\
 	& 
	=\sup_{s>0}\,\left(\int_s^\infty
	\tau^{-\frac{q}{r}}\Phi(\tau^{-1})^{q{\sigma}}\,d\tau\right)^{\frac{1}{q}}
	\left(\int_0^s \tau^{-\frac{q'}{r'}}\Phi(\tau^{-1})^{-q'{\sigma}}\,{d\tau}
	\right)^{\frac{1}{q'}}\\
	& 
	\le \sup_{s>0}\,\left\{C s^{\frac{1}{q}-\frac{1}{r}}\Phi(s^{-1})^\sigma
	\cdot Cs^{\frac{1}{q'}-\frac{1}{r'}}\Phi(s^{-1})^{-{\sigma}}\right\}<\infty.
	\end{aligned}
\end{equation}
Applying Lemma~\ref{Lemma:2.3} to \eqref{eq:3.8}, by \eqref{eq:3.9}
we obtain 
\begin{align*}
	\left\|S(t)\varphi\right\|_{L ^q(\log L)^\beta}^q
	 & \le C\sup_{s>0}\,\left\{s\Phi(s^{-1})^{\alpha}
	\varphi^{**}(s)^r\right\}^{\frac{q}{r}}
	\int_0^\infty\bigg(U_2(\tau)
	\int_0^\tau h_t^*(\xi)\,d\xi\bigg)^q\,d\tau\\
	 & \le 
	C\sup_{s>0}\,\left\{s\Phi(s^{-1})^{\alpha}
	\varphi^{**}(s)^r\right\}^{\frac{q}{r}}
  	\int_0^\infty\left(V_2(\tau)h_t^*(\tau)\right)^q\,d\tau\\
	 & 
	 =C\sup_{s>0}\,\left\{s\Phi(s^{-1})^{\alpha}
	\varphi^{**}(s)^r\right\}^{\frac{q}{r}}
  	\int_0^\infty\left(\tau^{1-\frac{1}{r}}\Phi(\tau^{-1})^{\sigma}
	h_t^*(\tau)\right)^q\,d\tau,\quad t>0.
\end{align*}
This together with Lemma~\ref{Lemma:3.2} implies that
$$
	\left\|S(t)\varphi\right\|_{L^q(\log L)^\beta}
	\le 
	Ct^{-\frac{N}{\theta}\left(\frac{1}{r}-\frac{1}{q}\right)}
	\Phi(t^{-1})^{\sigma}
	\sup_{s>0}\,\left\{s\Phi(s^{-1})^{\alpha}
	\varphi^{**}(s)^r\right\}^{\frac{1}{r}},
	\quad t>0.
$$
By Lemma~\ref{Lemma:3.3} and \eqref{eq:1.7}  
we see that Proposition~\ref{Proposition:3.1} holds in the case of $1\le r<q<\infty$. 
\vspace{3pt}
\newline
\underline{Step 2.}
Consider the case of $1\le r< q=\infty$. 
Let $\tilde{q}\in(r,\infty)$. 
By ($\Phi 2$) and \eqref{eq:2.8} 
we apply Proposition~\ref{Proposition:3.1} to obtain
\begin{align*}
	\|S(t)\varphi\|_{L^\infty} 
	& 
	=\left\|S\left(\frac{t}{2}\right)S\left(\frac{t}{2}\right)\varphi\right\|_{L^\infty}
	\\
 	& 
	\le Ct^{-\frac{N}{\theta \tilde{q}}}\left\|S\left(\frac{t}{2}\right)\varphi\right\|_{L^{\tilde{q}}}
 	=C^{-\frac{N}{\theta \tilde{q}}}\left\|S\left(\frac{t}{2}\right)\varphi\right\|_{L^{\tilde{q}}(\log  L)^0}
	\\
 	& 
	\le Ct^{-\frac{N}{\theta \tilde{q}}}\cdot Ct^{-\frac{N}{\theta}\left(\frac{1}{r}-\frac{1}{\tilde{q}}\right)}
	\Phi\biggr(\left(\frac{t}{2}\right)^{-1}\biggr)^{-\frac{\alpha}{r}}\Pi
	\le Ct^{-\frac{N}{\theta r}}\Phi(t^{-1})^{-\frac{\alpha}{r}}\Pi,
	\quad t>0,
\end{align*}
where
$$
	\Pi=
	\left\{
	\begin{array}{ll}
	\|\varphi\|_{{\mathfrak L}^{r,\infty}(\log{\mathfrak L})^\alpha}
	&\mbox{if $r\ge1$ and $\alpha\ge0$},
	\vspace{5pt}\\
	\|\varphi\|_{L^{r,\infty}(\log L)^\alpha}
	&\mbox{if $r>1$ and $\alpha\ge0$},
	\vspace{5pt}\\
	\|\varphi\chi_E\|_{L^{1,\infty}(\log L)^{\alpha+1}}
	&\mbox{if $\alpha>0$}.
	\end{array}
	\right.
$$
Thus Proposition~\ref{Proposition:3.1} follows in the case of $1\le r<q=\infty$. 
The proof of Proposition~\ref{Proposition:3.1} is complete.
$\Box$\vspace{5pt}

Then we apply the same arguments in the proof of \cite{IIK}*{Proposition~3.2} together with Proposition~\ref{Proposition:3.1}, 
and obtain the following proposition.
\begin{proposition}
\label{Proposition:3.2}
	Let $\theta\in(0,2]$, $1\le r< q\le \infty$, and $\alpha$, $\beta\ge 0$. 
	\begin{itemize}
	\item[\rm(1)]
		For any $r\ge1$, there exists $C_1>0$ such that
		\begin{equation*}
		\begin{split}
		 & \left|S(t)\varphi\right|_{q,\beta: T^{1/\theta}}
			+\|S(t)\varphi\|_{q,\beta;T^{1/\theta}}
			+|||S(t)\varphi|||_{q,\beta;T^{1/\theta}}\\
		 & \le 
			C_1t^{-\frac{N}{\theta}\left(\frac{1}{r}-\frac{1}{q}\right)}
			\Phi(t^{-1})^{-\frac{\alpha}{r}+\frac{\beta}{q}}
			|||\varphi|||_{r,\alpha;T^{1/\theta}}
		\end{split}
		\end{equation*}
		for $\varphi\in \mathfrak L^{r,\infty}_{{\rm ul}}(\log \mathfrak L)^\alpha$, $t\in(0,T]$, and $T\in(0,\infty]$.
	\item[\rm(2)]
		For any $r>1$, there exists $C_2>0$ such that
		\begin{equation*}
		\begin{split}
			 & \left|S(t)\varphi\right|_{q,\beta: T^{1/\theta}}
			+\|S(t)\varphi\|_{q,\beta;T^{1/\theta}}
			+|||S(t)\varphi|||_{q,\beta;T^{1/\theta}}\\
			 & \le 
			C_2t^{-\frac{N}{\theta}\left(\frac{1}{r}-\frac{1}{q}\right)}
			\Phi(t^{-1})^{-\frac{\alpha}{r}+\frac{\beta}{q}}
			\|\varphi\|_{r,\alpha;T^{1/\theta}}
		\end{split}			
		\end{equation*}
		for $\varphi\in L^{r,\infty}_{{\rm ul}}(\log L)^\alpha$, $t\in(0,T]$, and $T\in(0,\infty]$.
	\item[\rm(3)]
		Let $r=1$ and $\alpha>0$. 
		For any $T\in(0,\infty)$, there exists $C_3>0$ such that
		\begin{equation*}
		\begin{split}
			 & \left|S(t)\varphi\right|_{q,\beta: T^{1/\theta}}
			+\|S(t)\varphi\|_{q,\beta;T^{1/\theta}}
			+|||S(t)\varphi|||_{q,\beta;T^{1/\theta}}\\
			 & \le 
			C_3t^{-\frac{N}{\theta}\left(1-\frac{1}{q}\right)}
			\Phi(t^{-1})^{-\alpha+\frac{\beta}{q}}
			\|\varphi\|_{1,\alpha+1;T^{1/\theta}}
		\end{split}			
		\end{equation*}
		for $\varphi\in L^{1,\infty}_{{\rm ul}}(\log L)^{\alpha+1}$ and $t\in(0,T]$. 
	\end{itemize}
\end{proposition}

At the end of this section, 
we apply Hardy's inequality again to show that $L^q(\log L)^\alpha$ and $L^{q,\infty}(\log L)^\alpha$ are Banach spaces if $q>1$. 
\begin{lemma}
\label{Lemma:3.4}
Let $q\in(1,\infty)$ and $\alpha\in[0,\infty)$. 
For any $f\in{\mathcal L}$, 
set 
\begin{equation}
\label{eq:3.10}
\begin{split}
\|f\|_{L^q(\log L)^\alpha}' & :=\left(\int_0^\infty \Phi(s^{-1})^\alpha f^{**}(s)^q\,ds\right)^\frac{1}{q},\\
\|f\|_{L^{q,\infty}(\log L)^\alpha}' & :=
\sup_{s>0}\,\left\{s^{\frac{1}{q}}\Phi(s^{-1})^{\frac{\alpha}{q}}f^{**}(s)\right\}.
\end{split}			
\end{equation}
Then there exists $C>0$ such that 
\begin{gather}
\label{eq:3.11}
\|f\|_{L^q(\log L)^\alpha}\le \|f\|_{L^q(\log L)^\alpha}'\le C\|f\|_{L^q(\log L)^\alpha},\\
\label{eq:3.12}
\|f\|_{L^{q,\infty}(\log L)^\alpha}\le \|f\|_{L^{q,\infty}(\log L)^\alpha}'\le C\|f\|_{L^{q,\infty}(\log L)^\alpha},
\end{gather}
for $f\in{\mathcal L}$. 
Furthermore, 
$L^q(\log L)^\alpha$ and 
$L^{q,\infty}(\log L)^\alpha$ are Banach spaces equipped with norms $\|\cdot\|'_{L^q(\log L)^\alpha}$ and $\|\cdot\|'_{L^{q,\infty}(\log L)^\alpha}$, respectively.
\end{lemma}
{\bf Proof.}
Let $q\in(1,\infty)$ and $\alpha\in[0,\infty)$. 
It follows from (R4) and \eqref{eq:3.10} that 
\begin{align}
\label{eq:3.13}
 & \|f\|_{L^q(\log L)^\alpha}'
\ge\left(\int_0^\infty \Phi(s^{-1})^\alpha f^*(s)^q\,ds\right)^\frac{1}{q}=\|f\|_{L^q(\log L)^\alpha},\\
\label{eq:3.14}
 & \|f\|_{L^{q,\infty}(\log L)^\alpha}'
\ge\sup_{s>0}\left\{s^{\frac{1}{q}}\Phi(s^{-1})^{\frac{\alpha}{q}}f^*(s)\right\}=\|f\|_{L^{q,\infty}(\log L)^\alpha},
\end{align}
for $f\in{\mathcal L}$. 
Set $U_1(\tau)=\tau^{-1}\Phi(\tau^{-1})^{\frac{\alpha}{q}}$ 
and $V_1(\tau)=\Phi(\tau^{-1})^{\frac{\alpha}{q}}$ for $\tau>0$. 
If follows from Lemma~\ref{Lemma:3.1}~(1), (3) that
$$
\sup_{s>0}\left\{\|U_1\|_{L^q((s,\infty))}\|V_1^{-1}\|_{L^{q'}((0,s))}\right\}
\le \sup_{s>0}\left\{Cs^{-1+\frac{1}{q}}\Phi(s^{-1})^{\frac{\alpha}{q}}\cdot Cs^{\frac{1}{q'}}\Phi(s^{-1})^{-\frac{\alpha}{q}}\right\}<\infty.
$$
Then, by Lemma~\ref{Lemma:2.3} and \eqref{eq:3.10} we have
\begin{align*}
\|f\|_{L^q(\log L)^\alpha}'
 & =\left(\int_0^\infty  \left(U_1(s)\int_0^s f^*(\tau)\,d\tau\right)^q\,ds\right)^\frac{1}{q}\\
 & \le C\left(\int_0^\infty  \left(V_1(s)f^*(s)\right)^q\,ds\right)^\frac{1}{q}
=C\|f\|_{L^q(\log L)^\alpha},\quad f\in{\mathcal L}.
\end{align*}
This together with \eqref{eq:3.13} implies \eqref{eq:3.11}. 
Similarly, set 
$U_2(\tau):=\tau^{-1+\frac{1}{q}}\Phi(\tau^{-1})^{\frac{\alpha}{q}}$ and 
$V_2(\tau):=\tau^{\frac{1}{q}}\Phi(\tau^{-1})^{\frac{\alpha}{q}}$ for $\tau>0$. 
Since $U_2$ is non-increasing by $q>1$, it
follows from 
Lemma~\ref{Lemma:3.1}~(1) that 
$$
\sup_{s>0}\left\{\|U_2\|_{L^\infty((s,\infty))}\|V_2^{-1}\|_{L^1((0,s))}\right\}
\le \sup_{s>0}\left\{s^{-1+\frac{1}{q}}\Phi(s^{-1})^{\frac{\alpha}{q}}\cdot Cs^{1-\frac{1}{q}}\Phi(s^{-1})^{-\frac{\alpha}{q}}\right\}<\infty.
$$
Then Lemma~\ref{Lemma:2.3} implies that
$$
\|f\|_{L^{q,\infty}(\log L)^\alpha}'
=\sup_{s>0}\left\{U_2(s)\int_0^s f^*(\tau)\,d\tau\right\}\le C\sup_{s>0}\left\{V_2(s)f^*(s)\right\}=C\|f\|_{L^{q,\infty}(\log L)^\alpha}
$$
for $f\in{\mathcal L}$. This together with \eqref{eq:3.14} implies \eqref{eq:3.12}. 
Furthermore, we observe from (R5) that 
\begin{align*}
 & \|f\|_{L^q(\log L)^\alpha}'=\left(\int_0^\infty s^{-q}\Phi(s^{-1})^\alpha
 \left(\sup_{|E|=s}\int_E |f(x)|\,dx\right)^q\,ds\right)^\frac{1}{q},\\
 & \|f\|_{L^{q,\infty}(\log L)^\alpha}'=
\sup_{s>0}\sup_{|E|=s}\left\{s^{-1+\frac{1}{q}}\Phi(s^{-1})^{\frac{\alpha}{q}}\int_E |f(x)|\,dx\right\},
\quad f\in{\mathcal L}.
\end{align*}
Then we easily see that $\|\cdot\|_{L^q(\log L)^\alpha}'$ and $\|\cdot\|_{L^{q,\infty}(\log L)^\alpha}'$ are norms of $L^q(\log L)^\alpha$ and $L^{q,\infty}(\log L)^\alpha$, respectively. 
In addition, we observe
that $L^q(\log L)^\alpha$ and $L^{q,\infty}(\log L)^\alpha$ are the Banach spaces 
equipped with norms $\|\cdot\|_{L^q(\log L)^\alpha}'$. and $\|\cdot\|_{L^{q,\infty}(\log L)^\alpha}'$, respectively.  
Thus Lemma~\ref{Lemma:3.4} follows. 
$\Box$
\vspace{5pt}
\newline
By Lemma~\ref{Lemma:3.4} and \eqref{eq:1.9} we have:
\begin{lemma}
\label{Lemma:3.5}
Let $q\in(1,\infty)$ and $\alpha\in[0,\infty)$. 
For any $f\in{\mathcal L}$ and $\rho\in(0,\infty]$, set 
$$
|f|'_{q,\alpha;\rho}:=\sup_{z\in{\mathbb R}^N}\|f\chi_{B(z,\rho)}\|'_{L^q(\log L)^\alpha},
\quad
 \|f\|'_{q,\alpha;\rho}:=\sup_{z\in{\mathbb R}^N}\|f\chi_{B(z,\rho)}\|'_{L^{q,\infty}(\log L)^\alpha}.
$$
Then there exists $C>0$ such that 
$$
|f|_{q,\alpha;\rho}\le |f|'_{q,\alpha;\rho}\le C|f|_{q,\alpha;\rho},\quad
\|f\|_{q,\alpha;\rho}\le \|f\|'_{q,\alpha;\rho}\le C\|f\|_{q,\alpha;\rho},
\quad f\in{\mathcal L},\,\,\rho\in(0,\infty].
$$
Furthermore, 
$$
|f+g|'_{q,\alpha;\rho}\le |f|'_{q,\alpha;\rho}+|g|'_{q,\alpha;\rho},\quad
\|f+g\|'_{q,\alpha;\rho}\le \|f\|'_{q,\alpha;\rho}+\|g\|'_{q,\alpha;\rho},
\quad f, g\in{\mathcal L},\,\,\rho\in(0,\infty].
$$
\end{lemma}
\section{Proof of Theorem~\ref{Theorem:1.3}}
We introduce a definition of solutions to problem~\eqref{eq:P} in the sense of integral equations. 
\begin{definition}
\label{Definition:4.1}
Let $\theta\in(0,2]$, $p>1$, $T\in(0,\infty]$, and $\mu\in{\mathcal L}$. 
Set $F_p(s):=|s|^{p-1}s$ for $s\in{\mathbb R}$. 
Let $u$ be a measurable and finite almost everywhere function in ${\mathbb R}^N\times(0,T)$. 
We say that $u$ is a solution to problem~\eqref{eq:P} in ${\mathbb R}^N\times[0,T)$ if, 
for a.a.~$(x,t)\in{\mathbb R}^N\times(0,T)$, 
\begin{itemize}
  \item 
  $G(x-y,t-s)\mu(y)$ is integrable in ${\mathbb R}^N\times(0,t)$ with respect to $(y,s)\in{\mathbb R}^N\times(0,t)$;
  \item 
  $G(x-y,t-s)F_p(u(y,s))$ is integrable in ${\mathbb R}^N\times(0,t)$ with respect to $(y,s)\in{\mathbb R}^N\times(0,t)$;
  \item 
  $\displaystyle{
  u(x,t)=\int_0^t\int_{{\mathbb R}^N}G(x-y,t-s)\mu(y)\,dy\,ds+\int_0^t\int_{{\mathbb R}^N}G(x-y,t-s)F_p(u(y,s))\,dy\,ds.}$
\end{itemize}
\end{definition}
Similarly to \cite{HIT}*{Lemma~3.1}, we see that 
a solution to problem~\eqref{eq:P} in ${\mathbb R}^N\times(0,T)$ in the sense of Definition~\ref{Definition:4.1} is a solution in the sense of Definition~\ref{Definition:1.1}.
For the proofs of Theorems~\ref{Theorem:1.3} and~\ref{Theorem:1.4}, 
we find a solution in the sense of  Definition~\ref{Definition:4.1}. 

In this section we apply the contraction mapping theorem to problem~\eqref{eq:P} in uniformly local weak Zygmund spaces, 
and prove Theorem~\ref{Theorem:1.3}. We also prove Corollaries~\ref{Corollary:1.1} and \ref{Corollary:1.2}. 
One of the main ingredients of the proof of Theorem~\ref{Theorem:1.3} is the following integral estimate, 
which is proved by the real interpolation method in uniformly local Zygmund spaces.  
\begin{Proposition}
\label{Proposition:4.1}
Let $\theta\in(0,2]$ with $N>\theta$. 
Then there exists $C_*>0$ such that
\begin{equation}
\label{eq:4.1}
	\sup_{0<t<T}\bigg\|\int_0^t S(s)\mu\,ds\,\bigg\|'_{N/(N-\theta),N/\theta;T^{1/\theta}}
	\le C_*|||\mu|||_{1,(N-\theta)/\theta;T^{1/\theta}}
\end{equation}
for $\mu\in {\mathfrak L}^{1,\infty}_{{\rm ul}}(\log{\mathfrak L})^{(N-\theta)/\theta}$ and $T\in(0,\infty]$.
\end{Proposition}
For the proof of Proposition~\ref{Proposition:4.1}, 
we introduce some notations.
Let $X_0$ and $X_1$ be Banach spaces.
For any $f\in X_0+X_1$ and $\lambda>0$, we define
\begin{equation}
\label{eq:4.2}
	K(f,\lambda; X_0, X_1)
	:=\inf\{\|f_0\|_{X_0}+\lambda\|f_1\|_{X_1}\,:\,f=f_0+f_1,\,\, f_0\in X_0,\,\, f_1\in X_1\}.
\end{equation}
For any $\kappa\in(0,1)$, we define 
$$
	(X_0,X_1)_{\kappa,\infty}
	:=\{f\in X_0+X_1\,:\, \|f\|_{(X_0,X_1)_{\kappa,\infty}}<\infty\},
$$
where
\begin{equation}
\label{eq:4.3}
	\|f\|_{(X_0,X_1)_{\kappa,\infty}}
	:=\sup_{\lambda>0}\,\lambda^{-\kappa}K(f,\lambda; X_0, X_1)
\end{equation}
(see e.g., \cite{BS}*{Chaper~5, Section~1}).
Then we have:
\begin{Lemma}
\label{Lemma:4.1}
Let
\begin{equation}
\label{eq:4.4}
	1\le q_0<q<q_1<\infty,
	\quad
	\frac{1}{q}=\frac{1-\kappa}{q_0}+\frac{\kappa}{q_1}
	\,\,\mbox{ with }\,\,\kappa\in(0,1),
	\quad
	\alpha\in[0,\infty).
\end{equation}
Set 
$$
X_0={\mathfrak L}^{q_0,\infty}(\log{\mathfrak L})^{\alpha \frac{q_0}{q}},\quad
X_1={\mathfrak L}^{q_1,\infty}(\log{\mathfrak L})^{\alpha \frac{q_1}{q}}.
$$
Then
\begin{equation}
\label{eq:4.5}
	\|f\|_{L^{q,\infty}(\log L)^\alpha}
	\le 2^{\frac{1}{q_0}}\|f\|_{(X_0,X_1)_{\kappa,\infty}},
	\quad 
	f\in (X_0,X_1)_{\kappa,\infty}.
\end{equation}
\end{Lemma}
\begin{remark}
\label{Remark:4.1}
In the case $\alpha=0$,
it follows from \eqref{eq:1.8} that
\[
	({\mathfrak L}^{q_0,\infty}(\log{\mathfrak L})^0,
	{\mathfrak L}^{q_1,\infty}(\log{\mathfrak L})^0)_{\kappa,\infty}
	=(L^{q_0}, L^{q_1})_{\kappa, \infty} = L^{q,\infty}
\]
{\rm ({\it see e.g.,} \cite{BL}*{{\it Theorem} 5.2.1})}. 
This relation is used in the proof of Theorem~{\rm\ref{Theorem:1.4}} {\rm({\it see Section~$5$})}.
\end{remark}
{\bf Proof.}
Assume \eqref{eq:4.4}. Let 
\begin{align*}
 & f\in\left({\mathfrak L}^{q_0,\infty}(\log{\mathfrak L})^{\alpha\frac{q_0}{q}},{\mathfrak L}^{q_1,\infty}(\log{\mathfrak L})^{\alpha\frac{q_1}{q}}\right)_{\kappa,\infty},\\
 & f=f_0+f_1\in {\mathfrak L}^{q_0,\infty}(\log{\mathfrak L})^{\alpha\frac{q_0}{q}}+{\mathfrak L}^{q_1,\infty}(\log{\mathfrak L})^{\alpha\frac{q_1}{q}}.
\end{align*}
It follows from Jensen's inequality, (R2), and (R4) that
\begin{equation}
\label{eq:4.6}
\begin{split}
	f^*(s)
	&
	\le f^*_0\left(\frac{s}{2}\right)+f^*_1\left(\frac{s}{2}\right)
	\le \frac{2}{s}\int_0^{s/2}f^*_0(\tau)\,d\tau+\frac{2}{s}\int_0^{s/2}f^*_1(s)\,d\tau
	\\
	&
	\le\bigg(\frac{2}{s}\int_0^{s/2}f^*_0(\tau)^{q_0}\,d\tau\bigg)^{\frac{1}{q_0}}
	+\bigg(\frac{2}{s}\int_0^{s/2}f^*_1(\tau)^{q_1}\,d\tau\bigg)^{\frac{1}{q_1}}
	\\
	 & \le 2^{\frac{1}{q_0}}\bigg\{s^{-\frac{1}{q_0}}\bigg(\int_0^s f^*_0(\tau)^{q_0}\,d\tau\bigg)^{\frac{1}{q_0}}
	+s^{-\frac{1}{q_1}}\bigg(\int_0^s f^*_1(\tau)^{q_1}\,d\tau\bigg)^{\frac{1}{q_1}}\bigg\},
	\quad s>0.
\end{split}
\end{equation}
Since $\Phi(s^{-1})\ge 1$ for $s\in(0,\infty)$ (see ($\Phi$1)), 
by \eqref{eq:2.1} and \eqref{eq:4.6} we have
$$
	\begin{aligned}
	&
	\left\{s\Phi(s^{-1})^\alpha f^*(s)^q\right\}^\frac{1}{q}
	=s^{\frac{1-\kappa}{q_0}+\frac{\kappa}{q_1}}\Phi(s^{-1})^{\frac{\alpha}{q}} f^*(s)
	\\
	&
	\le  2^{\frac{1}{q_0}}\bigg\{
	s^{-\kappa\left(\frac{1}{q_0}-\frac{1}{q_1}\right)}
	\Phi(s^{-1})^{\frac{\alpha}{q}}
	\bigg(\int_0^s f^*_0(\tau)^{q_0}\,d\tau\bigg)^{\frac{1}{q_0}}
	\\
	&\hspace{4cm}
	+s^{(1-\kappa)\left(\frac{1}{q_0}-\frac{1}{q_1}\right)}\Phi(s^{-1})^{\frac{\alpha}{q}}
	\bigg(\int_0^s f^*_1(\tau)^{q_1}\,d\tau\bigg)^{\frac{1}{q_1}}\bigg\}
	\\
	&
	= 2^{\frac{1}{q_0}}s^{-\kappa\left(\frac{1}{q_0}-\frac{1}{q_1}\right)}
	\bigg\{\bigg(\Phi(s^{-1})^{\frac{\alpha q_0}{q}}
	\int_0^s f^*_0(\tau)^{q_0}\,d\tau\bigg)^{\frac{1}{q_0}}
	\\
	&\hspace{4cm}
	+s^{\frac{1}{q_0}-\frac{1}{q_1}}\bigg(\Phi(s^{-1})^{\frac{\alpha q_1}{q}}
	\int_0^s f^*_1(\tau)^{q_1}\,d\tau\bigg)^{\frac{1}{q_1}}\bigg\}
	\\
	&
	\le 2^{\frac{1}{q_0}}s^{-\kappa\left(\frac{1}{q_0}-\frac{1}{q_1}\right)}
	\bigg(\|f_0\|_{{\mathfrak L}^{q_0,\infty}(\log{\mathfrak L})^{\alpha \frac{q_0}{q}}}
	+s^{\frac{1}{q_0}-\frac{1}{q_1}}\|f_1\|_{{\mathfrak L}^{q_1,\infty}(\log{\mathfrak L})^{\alpha \frac{q_1}{q}}}\bigg),
	\quad s>0. 
	\end{aligned}
$$
Then it follows from \eqref{eq:4.2} and \eqref{eq:4.3} that
$$
	\begin{aligned}
	&
	\left\{s\Phi(s^{-1})^\alpha f^*(s)^q\right\}^{\frac{1}{q}}
	\\
	&
	\le 2^{\frac{1}{q_0}}s^{-\kappa\left(\frac{1}{q_0}-\frac{1}{q_1}\right)} 
	K\bigg(f, s^{\frac{1}{q_0}-\frac{1}{q_1}};{\mathfrak L}^{q_0,\infty}(\log{\mathfrak L})^{\alpha \frac{q_0}{q}},
	{\mathfrak L}^{q_1,\infty}(\log{\mathfrak L})^{\alpha \frac{q_1}{q}} \bigg)
	\\
	&
	\le 2^{\frac{1}{q_0}}\sup_{\lambda>0}\lambda^{-\kappa} 
	K\bigg(f, \lambda;{\mathfrak L}^{q_0,\infty}(\log{\mathfrak L})^{\alpha \frac{q_0}{q}},
	{\mathfrak L}^{q_1,\infty}(\log{\mathfrak L})^{\alpha \frac{q_1}{q}} \bigg)
	=2^{\frac{1}{q_0}}\|f\|_{(X_0,X_1)_{\kappa,\infty}},
	\quad s>0.
	\end{aligned}
$$
This together with \eqref{eq:1.5} implies \eqref{eq:4.5}.
Thus Lemma~\ref{Lemma:4.1} follows.
$\Box$
\vspace{5pt}

\noindent
{\bf Proof of Proposition~\ref{Proposition:4.1}.}
We employ the arguments in the proofs of \cite{Meyer}*{Theorem 18.1} and \cite{OT}*{Lemma 4.1}. 
Let $z\in{\mathbb R}^N$ and $t\in(0,T)$. 
Set 
\begin{equation}
\label{eq:4.7}
\begin{split}
f_z(x,t) & :=\chi_{B(z,T^{1/\theta})}(x)\int_0^t [S(s)\mu](x)\,ds\\
 & \, =\chi_{B(z,T^{1/\theta})}(x)\int_0^\infty \chi_{(0,t)}(s)[S(s)\mu](x)\,ds,\quad x\in{\mathbb R}^N.
\end{split}
\end{equation}
For any $\eta>0$, we set 
\begin{equation}
\label{eq:4.8}
\begin{split}
 & f_{z,0}^\eta(x,t):=\chi_{B(z,T^{1/\theta})}(x)\int_0^\eta \chi_{(0,t)}(s)[S(s)\mu](x)\,ds,\\
 & f_{z,1}^\eta(x,t):=\chi_{B(z,T^{1/\theta})}(x)\int_\eta^\infty \chi_{(0,t)}(s)[S(s)\mu](x)\,ds,\quad (x,t)\in{\mathbb R}^N\times(0,T).
\end{split}
\end{equation}
Then 
\begin{equation}
\label{eq:4.9}
f_z(x,t)=f_{z,0}^\eta(x,t)+f_{z,1}^\eta(x,t),\quad x\in{\mathbb R}^N.
\end{equation}
Let $q_0$, $q_1\in(1,\infty)$ be such that
\begin{equation}
\label{eq:4.10}
	1<q_0<p_*=\frac{N}{N-\theta}<q_1<\infty,\qquad
	\frac{1}{p_*}=\frac{1/2}{q_0}+\frac{1/2}{q_1},
\end{equation}
and set 
\begin{equation}
\label{eq:4.11}
X_0=\mathfrak L^{q_0,\infty}(\log \mathfrak L)^{q_0\gamma},
\quad 
X_1=\mathfrak L^{q_1,\infty}(\log \mathfrak L)^{q_1\gamma},
\quad 
\gamma=\frac{N-\theta}{\theta}.
\end{equation}
Since
\begin{equation}
\label{eq:4.12}
	\frac{N}{\theta}\bigg(1-\frac{1}{q_0}\bigg)<\frac{N}{\theta}\bigg(1-\frac{1}{p_*}\bigg)=1<\frac{N}{\theta}\bigg(1-\frac{1}{q_1}\bigg),
\end{equation}
by Proposition~\ref{Proposition:3.2}~(1), \eqref{eq:1.9}, and \eqref{eq:4.8} we have
\begin{align*}
\|f_{z,0}^\eta(t)\|_{X_0}
 & \le\bigg|\bigg|\bigg|\,\int_0^\eta \chi_{(0,t)}(s)S(s)\mu\,ds\,\bigg|\bigg|\bigg|_{q_0,q_0\gamma;T^{1/\theta}}\\
 & \le\int_0^\eta \chi_{(0,t)}(s)|||S(s)\mu|||_{q_0,q_0\gamma;T^{1/\theta}}\,ds\\
 & \le C\int_0^\eta s^{-\frac{N}{\theta}\left(1-\frac{1}{q_0}\right)}|||\mu|||_{1,\gamma;T^{1/\theta}}\,ds
 \le C\eta^{1-\frac{N}{\theta}\left(1-\frac{1}{q_0}\right)}|||\mu|||_{1,\gamma;T^{1/\theta}},\\
\|f_{z,1}^\eta(t)\|_{X_1}
 & \le\bigg|\bigg|\bigg|\,\int_\eta^\infty \chi_{(0,t)}(s)S(s)\mu\,ds\,\bigg|\bigg|\bigg|_{q_1,q_1\gamma;T^{1/\theta}}\\
 & \le\int_\eta^\infty \chi_{(0,t)}(s)|||S(s)\mu|||_{q_1,q_1\gamma;T^{1/\theta}}\,ds\\
 & \le C\int_\eta^\infty s^{-\frac{N}{\theta}\left(1-\frac{1}{q_1}\right)}|||\mu|||_{1,\gamma;T^{1/\theta}}\,ds
 \le C\eta^{1-\frac{N}{\theta}\left(1-\frac{1}{q_1}\right)}|||\mu|||_{1,\gamma;T^{1/\theta}}.
\end{align*}
These together with \eqref{eq:4.2} and \eqref{eq:4.9} imply that
\begin{equation}
\label{eq:4.13}
\begin{split}
K(f_z(t),\lambda; X_0,X_1)
 & \le\|f_{z,0}^\eta\|_{X_0}+\lambda\|f_{z,1}^\eta\|_{X_1}\\
 & \le C\left(\eta^{1-\frac{N}{\theta}\left(1-\frac{1}{q_0}\right)}+\lambda\eta^{1-\frac{N}{\theta}\left(1-\frac{1}{q_1}\right)}\right)|||\mu|||_{1,\gamma;T^{1/\theta}}
\end{split}
\end{equation}
for $z\in{\mathbb R}^N$, $t\in(0,T)$, $\lambda>0$, and $\eta>0$.
Taking $\eta>0$ so that 
\begin{equation}
\label{eq:4.14}
\eta^{1-\frac{N}{\theta}\left(1-\frac{1}{q_0}\right)}=\lambda^\frac{1}{2},
\end{equation}
since $p_*=N/(N-\theta)$, by \eqref{eq:4.10} we have 
\begin{equation}
\label{eq:4.15}
	\eta^{1-\frac{N}{\theta}\left(1-\frac{1}{q_1}\right)}
	=\eta^{1-\frac{N}{\theta}\left(1-\frac{2}{p_*}+\frac{1}{q_0}\right)}
	=\eta^{-1+\frac{N}{\theta}\left(1-\frac{1}{q_0}\right)}
	=\lambda^{-\frac{1}{2}}.
\end{equation}
Combing \eqref{eq:4.13}, \eqref{eq:4.14}, and \eqref{eq:4.15}, we obtain 
$$
K(f_z(t),\lambda; X_0,X_1)\le C\lambda^{\frac{1}{2}}|||\mu|||_{1,\gamma;T^{1/\theta}},\quad \lambda>0,
$$
for $z\in{\mathbb R}^N$ and $t\in(0,T)$. 
Then, by \eqref{eq:4.3} and \eqref{eq:4.10} we apply Lemma~\ref{Lemma:4.1} with $\kappa=1/2$ to obtain 
\begin{align*}
 & \bigg\|\int_0^t S(s)\mu\,ds\,\bigg\|_{N/(N-\theta),N/\theta;T^{1/\theta}}\\
 & =\sup_{z\in{\mathbb R}^N}\|f_z\|_{L^{p_*,\infty}(\log L)^{p_*\gamma}}
\le C\sup_{z\in{\mathbb R}^N}\|f_z\|_{(X_0,X_1)_{1/2,\infty}}\\
 & =C\sup_{z\in{\mathbb R}^N}\sup_{\lambda>0}\left\{\lambda^{-\frac{1}{2}}K(f_z,\lambda; X_0,X_1)\right\}
\le C|||\mu|||_{1,\gamma;T^{1/\theta}},\quad t\in(0,T). 
\end{align*}
This together with Lemma~\ref{Lemma:3.5} implies \eqref{eq:4.1}, and Proposition~\ref{Proposition:4.1} follows.
$\Box$
\vspace{8pt}

We are ready to prove Theorem~\ref{Theorem:1.3}. 
For the proof, let $\theta<N$ and set  
\begin{equation}
\label{eq:4.16}
T_*\in(0,\infty),\quad 
p=p_* := \frac{N}{N-\theta},\quad \gamma:=\frac{N}{\theta p_*}=\frac{N-\theta}{\theta},
\quad\mu\in {\mathfrak L}^{1,\infty}_{{\rm ul}}(\log{\mathfrak L})^\gamma.
\end{equation}
Define
$$
	X_T:=C((0,T): L_{{\rm ul}}^{p,\infty}(\log L)^{p\gamma}),\quad T\in(0,T_*).
$$
Let $C^*=2C_*$, where $C_*$ is given in Proposition~\ref{Proposition:4.1}.
For 
$\varepsilon>0$, 
we define
\begin{equation}
\label{eq:4.17}
X_T(C^*\varepsilon):=
\left\{u\in X_T\,:\,\sup_{0<t<T}\,\|u(t)\|'_{p,p\gamma;T^{1/\theta}}\le C^*\varepsilon\right\}.
\end{equation}
For any $u$, $v\in X_T$, set 
$$
	d_X(u,v) :=\sup_{0<t<T}\,\|u(t)-v(t)\|'_{p,p\gamma;T^{1/\theta}}.
$$
Then $(X_T,d_X)$ is a complete metric space and $X_T(C^*\varepsilon)$ is closed in $(X_T,d_X)$. 
Define
\begin{equation}
\label{eq:4.18}
	{\mathcal F}(u):=\int_0^t S(t-s)\mu\,ds+\int_0^t S(t-s)F_p(u(s))\,ds,\quad u\in X_T(C^*\varepsilon),
\end{equation}
where $F_p(\tau)=|\tau|^{p-1}\tau$ for $\tau\in \mathbb{R}$.
\begin{lemma}
\label{Lemma:4.2}
Let $\theta\in(0,2]$ with $N>\theta$. 
Let $p$, $T_*$, and $\gamma$ be as in \eqref{eq:4.16}, and $T\in(0,T_*)$. 
Then there exists $C>0$ such that 
$$
d_X({\mathcal F}(u),{\mathcal F}(v))\le C\varepsilon^{p-1}\,d_X(u,v)
\quad\mbox{for}\quad u,v\in X_T(C^*\varepsilon),\,\,\varepsilon>0.
$$
\end{lemma}
{\bf Proof.}
Let $u$, $v\in X_T(C^*\varepsilon)$. Let $0<s<t<T$. 
It follows that
\begin{equation}
\label{eq:4.19}
	|F_p(u(x,s))-F_p(v(x,s))|\le w(x,s)|u(x,s)-v(x,s)|,\quad x\in{\mathbb R}^N, 
\end{equation}
where $w(x,s):=p(|u(x,s)|^{p-1}+|v(x,s)|^{p-1})$.
By Lemmas~\ref{Lemma:2.1}, \ref{Lemma:2.2}, and \ref{Lemma:3.5} we have 
$$
	\begin{aligned}
 	& 
	\|F_p(u(s))-F_p(v(s))\|_{1,p\gamma;T^{1/\theta}}
	\\
 	& 
	\le C\|w(s)\|_{\frac{p}{p-1},p\gamma;T^{1/\theta}} \|u(s)-v(s)\|_{p,p\gamma;T^{1/\theta}}
	\\
	& 
	\le C\left(\|u(s)\|_{p,p\gamma;T^{1/\theta}}^{p-1}+\|v(s)\|_{p,p\gamma;T^{1/\theta}}^{p-1}\right)
	\|u(s)-v(s)\|_{p,p\gamma;T^{1/\theta}}.
	\end{aligned}
$$
This together with Lemma~\ref{Lemma:3.5} and \eqref{eq:4.17} implies that
\begin{equation}
\label{eq:4.20}
	\sup_{0<s<t}\,\|F_p(u(s))-F_p(v(s))\|_{1,p\gamma;T^{1/\theta}}
	\le C(C^*\varepsilon)^{p-1}d_X(u,v).
\end{equation}

Let $z\in{\mathbb R}^N$, and set 
\begin{equation}
\label{eq:4.21}
\begin{split}
g_z(x,t) & :=\chi_{B(z,T^{1/\theta})}\int_0^tS(t-s)[F_p(u(s))-F_p(v(s))]\,ds\\
 & \,=\chi_{B(z,T^{1/\theta})}\int_0^\infty \chi_{(0,t)}(s)S(t-s)[F_p(u(s))-F_p(v(s))]\,ds\\
 & \, =\chi_{B(z,T^{1/\theta})}\int_0^\infty \chi_{(0,t)}(s)S(s)[F_p(u(t-s))-F_p(v(t-s))]\,ds
\end{split}
\end{equation}
for $(x,t)\in{\mathbb R}^N\times(0,T)$. For any $\eta>0$, we set 
\begin{equation}
\label{eq:4.22}
\begin{split}
 & g_{z,0}^\eta(x,t):=\chi_{B(z,T^{1/\theta})}\int_0^\eta \chi_{(0,t)}(s)S(s)[F_p(u(t-s))-F_p(v(t-s))]\,ds,\\
 & g_{z,1}^\eta(x,t):=\chi_{B(z,T^{1/\theta})}\int_\eta^\infty \chi_{(0,t)}(s)S(s)[F_p(u(t-s))-F_p(v(t-s))]\,ds, 
\end{split}
\end{equation}
for $(x,t)\in{\mathbb R}^N\times(0,T)$. 
Then 
\begin{equation}
\label{eq:4.23}
g_z(x,t)=g_{z,0}^\eta(x,t)+g_{z,1}^\eta(x,t),\quad (x,t)\in{\mathbb R}^N\times(0,T). 
\end{equation}
Let $q_0$, $q_1\in(0,\infty)$ (resp.~$X_0$ and $X_1$) be as \eqref{eq:4.10} (resp.~\eqref{eq:4.11}).
Since $p\gamma=\gamma+1$ (see \eqref{eq:4.16}),
by \eqref{eq:4.12} and \eqref{eq:4.20} 
we apply Proposition~\ref{Proposition:3.2}~(3) to obtain 
\begin{align*}
\|g_{z,0}^\eta(t)\|_{X_0} 
 & \le\bigg|\bigg|\bigg|\int_0^\eta \chi_{(0,t)}(s)\,S(s)\{F_p(u(t-s))-F_p(v(t-s))\}\,ds\,\bigg|\bigg|\bigg|_{q_0,q_0\gamma;T^{1/\theta}}\\
 & \le\int_0^\eta \chi_{(0,t)}(s)\,|||S(s)\{F_p(u(t-s))-F_p(v(t-s))\}|||_{q_0,q_0\gamma;T^{1/\theta}}\,ds\\
 & \le C\int_0^\eta \chi_{(0,t)}(s)\, s^{-\frac{N}{\theta}\left(1-\frac{1}{q_0}\right)}\|F_p(u(t-s))-F_p(v(t-s))\|_{1,\gamma+1;T^{1/\theta}}\,ds\\
 & \le C\int_0^\eta s^{-\frac{N}{\theta}\left(1-\frac{1}{q_0}\right)}\,ds\,\cdot \sup_{0<s<t}\,\|F_p(u(s))-F_p(v(s))\|_{1,p\gamma;T^{1/\theta}}\\
 & \le C(C^*\varepsilon)^{p-1}\eta^{1-\frac{N}{\theta}\left(1-\frac{1}{q_0}\right)}d_X(u,v),\\
\|g_{z,1}^\eta(t)\|_{X_1}
 & \le\bigg|\bigg|\bigg|\int_\eta^\infty \chi_{(0,t)}(s)\,S(s)\{F_p(u(t-s))-F_p(v(t-s))\}\,ds\,\bigg|\bigg|\bigg|_{q_1,q_1\gamma;T^{1/\theta}}\\
 & \le\int_\eta^\infty\chi_{(0,t)}(s)\,|||S(s)\{F_p(u(t-s))-F_p(v(t-s))\}|||_{q_1,q_1\gamma;T^{1/\theta}}\,ds\\
 & \le C\int_\eta^\infty\chi_{(0,t)}(s)\,s^{-\frac{N}{\theta}\left(1-\frac{1}{q_1}\right)}\|F_p(u(t-s))-F_p(v(t-s))\|_{1,\gamma+1;T^{1/\theta}}\,ds\\
 & \le C\int_\eta^\infty s^{-\frac{N}{\theta}\left(1-\frac{1}{q_1}\right)}\,ds\,
 \cdot
 \sup_{0<s<t}\,\|F_p(u(s))-F_p(v(s))\|_{1,p\gamma;T^{1/\theta}}\\
 & \le C(C^*\varepsilon)^{p-1}\eta^{1-\frac{N}{\theta}\left(1-\frac{1}{q_1}\right)}d_X(u,v).
\end{align*}
These together with \eqref{eq:4.2} and \eqref{eq:4.23} imply that
\begin{equation}
\label{eq:4.24}
\begin{split}
K(g_z(t),\lambda; X_0,X_1)
 & \le\|g_{z,0}^\eta(t)\|_{X_0}+\lambda\|g_{z,1}^\eta(t)\|_{X_1}\\
 & \le C(C^*\varepsilon)^{p-1}\left(\eta^{1-\frac{N}{\theta}\left(1-\frac{1}{q_0}\right)}
 +\lambda\eta^{1-\frac{N}{\theta}\left(1-\frac{1}{q_1}\right)}\right)d_X(u,v)
\end{split}
\end{equation}
for $z\in{\mathbb R}^N$, $t\in(0,T)$, $\lambda>0$, and $\eta>0$.
Taking $\eta>0$ satisfying \eqref{eq:4.14}, 
by \eqref{eq:4.15} and \eqref{eq:4.24} we have 
$$
K(g_z(t),\lambda; X_0,X_1)\le C(C^*\varepsilon)^{p-1}\lambda^{\frac{1}{2}}d_X(u,v)
$$
for $z\in{\mathbb R}^N$, $t\in(0,T)$, and $\lambda>0$. 
Then, by \eqref{eq:4.3} and \eqref{eq:4.10} we apply Lemma~\ref{Lemma:3.4} and Lemma~\ref{Lemma:4.1} with $\kappa=1/2$ to obtain 
\begin{align*}
 & d_X({\mathcal F}(u),{\mathcal F}(v))\\
 & =\sup_{t\in(0,T)}\bigg\|\int_0^tS(t-s)\{F_p(u(s))-F_p(v(s))\}\,ds\,\bigg\|'_{p,p\gamma;T^{1/\theta}}\\
 & \le C\sup_{t\in(0,T)}\sup_{z\in{\mathbb R}^N}\|g_z(t)\|_{L^{p,\infty}(\log L)^{p\gamma}}
\le C\sup_{t\in(0,T)}\sup_{z\in{\mathbb R}^N} \|g_z(t)\|_{(X_0,X_1)_{1/2,\infty}}\\
 & =C\sup_{t\in(0,T)}\sup_{z\in{\mathbb R}^N}\sup_{\lambda>0}\left\{\lambda^{-\frac{1}{2}}K(g_z(t),\lambda; X_0,X_1)\right\}
 \le C(C^*\varepsilon)^{p-1}d_X(u,v).
\end{align*}
Thus we obtain the desired inequality. 
The proof is complete. 
$\Box$
\vspace{5pt}

\noindent
{\bf Proof of Theorem~\ref{Theorem:1.3}.}
Let $p$, $T_*$, and $\gamma$ be as in \eqref{eq:4.16}. 
Let $\varepsilon>0$ be small enough,
and let $\mu\in {\mathfrak L}^{1,\infty}_{{\rm ul}}(\log{\mathfrak L})^\gamma$ be such that 
$|||\mu|||_{1,\gamma;T^{1/\theta}}<\varepsilon$ for some $T\in(0,T_*]$. 
By Proposition~\ref{Proposition:4.1} and Lemma~\ref{Lemma:4.2} we have
$$
	\begin{aligned}
	\sup_{t\in(0,T)}\|{\mathcal F}(u(t))\|'_{p,p\gamma;T^{1/\theta}}
	& 
	\le \sup_{0<t<T}\bigg\|\int_0^tS(s)\mu\,ds\bigg\|'_{p,p\gamma,T^{1/\theta}}
	+{d_X\bigl({\mathcal F}(u),{\mathcal F}(0)\bigr)}
	\\
	& 
	\le C_*\varepsilon+C\varepsilon^{p-1}d_X(u,0)
	 \le C_*\varepsilon+C\varepsilon^{p-1}\cdot C^*\varepsilon\le C^*\varepsilon
	\end{aligned}
$$
for $u\in X_T(C^*\varepsilon)$. 
This implies that ${\mathcal F}(u)\in X_T(C^*\varepsilon)$ for $u\in X_T(C^*\varepsilon)$.
Furthermore, 
taking small enough $\varepsilon>0$ if necessary, 
by Lemma~\ref{Lemma:4.2} we have 
$$
	d_X({\mathcal F}(u),{\mathcal F}(v))
	\le C\varepsilon^{p-1}\,d_X(u,v)\le\frac{1}{2}d_X(u,v)
$$
for $u$, $v\in X_T(C^*\varepsilon)$. 
We apply the contraction mapping theorem to find a unique $u_*\in X_T(C^*\varepsilon)$ such that 
${\mathcal F}(u_*)=u_*$ in $X_T(C^*\varepsilon)$. 
Then $u_*$ is a solution to problem~\eqref{eq:P} in ${\mathbb R}^n\times[0,T)$ 
in the sense of Definition~\ref{Definition:4.1}, with $u_*$ satisfying \eqref{eq:1.11}. 
Thus Theorem~\ref{Theorem:1.3} follows.~$\Box$
\vspace{5pt}
\newline
{\bf Proof of Corollary~\ref{Corollary:1.1}}
Since $\alpha>\gamma$, 
it follows from \eqref{eq:2.3} that 
$$
	|||\mu|||_{1,\gamma;T^{1/\theta}}
	\le C\Phi(T^{-1/\theta})^{\gamma-\alpha}
	|||\mu|||_{1,\alpha;T^{1/\theta}}\to 0
	\quad\mbox{as}\quad T\to +0.
$$
Then, by Theorem~\ref{Theorem:1.3} we find a solution $u$ to problem~\eqref{eq:P} in ${\mathbb R}^N\times(0,T)$ 
for small enough $T>0$.  
Thus Corollary~\ref{Corollary:1.1} follows.
$\Box$\vspace{5pt}
\newline
{\bf Proof of Corollary~\ref{Corollary:1.2}.} 
It follows from \eqref{eq:1.12} that $\mu_c\in {\mathfrak L}^{1,\infty}_{{\rm ul}}(\log{\mathfrak L})^{\gamma}$. 
Then Corollary~\ref{Corollary:1.2} follows from Theorem~\ref{Theorem:1.3}.
$\Box$ 
\section{Proof of Theorem~\ref{Theorem:1.4}}
We modify the arguments in Section~4 to prove Theorem~\ref{Theorem:1.4}. 
We prepare the following proposition, instead of Proposition~\ref{Proposition:4.1}. 
\begin{Proposition}
\label{Proposition:5.1}
Let $\theta\in(0,2]$ with $N>\theta$. 
Let 
$$
p>p_* := \frac{N}{N-\theta},\quad r=r_*:=\frac{N(p-1)}{\theta p}.
$$
There exists $C_{**}>0$ such that
\begin{equation*}
	\sup_{0<t<T}\bigg\|\int_0^t S(s)\mu\,ds\,\bigg\|'_{pr,0;T^{1/\theta}}
	\le C_{**}\|\mu\|_{r,0;T^{1/\theta}}
\end{equation*}
for $\mu\in L^{r,\infty}_{{\rm ul}}$ and $T\in(0,\infty]$. 
\end{Proposition}
{\bf Proof.}
For any $z\in{\mathbb R}^N$ and $\eta>0$, 
let $f_z$, $f_{z,0}^\eta$, and $f_{z,1}^\eta$ be as in \eqref{eq:4.7} and \eqref{eq:4.8}. 
Let $q_0$, $q_1\in(1,\infty)$ be such that 
\begin{equation}
\label{eq:5.1}
1<r<q_0<pr<q_1<\infty,\qquad
\frac{1}{pr}=\frac{1/2}{q_0}+\frac{1/2}{q_1}.
\end{equation}
Since
\begin{equation}
\label{eq:5.2}
\frac{N}{\theta}\bigg(\frac{1}{r}-\frac{1}{pr}\bigg)=\frac{N}{\theta r}\frac{p-1}{p}=1,
\quad
\frac{N}{\theta}\bigg(\frac{1}{r}-\frac{1}{q_0}\bigg)<1<\frac{N}{\theta}\bigg(\frac{1}{r}-\frac{1}{q_1}\bigg),
\end{equation}
by Proposition~\ref{Proposition:3.2}~(2) and \eqref{eq:1.9} we see that
\begin{align*}
	\|f_{z,0}^\eta\|_{L^{q_0}}
	&
	\le
	\int_0^\eta \chi_{(0,t)}(s)|S(s)\mu|_{q_0,0;T^{1/\theta}}\,ds
	\\
	&
	\le C \int_0^\eta s^{-\frac{N}{\theta}\left(\frac{1}{r}-\frac{1}{q_0}\right)}\|\mu\|_{r,0;T^{1/\theta}}\,ds
	\le C\eta^{1-\frac{N}{\theta}\left(\frac{1}{r}-\frac{1}{q_0}\right)}\,
	\|\mu\|_{r,0;T^{1/\theta}},\\
	\|f_{z,1}^\eta\|_{L^{q_1}}
	&
	\le
	\int_\eta^\infty \chi_{(0,t)}(s)|S(s)\mu|_{q_1,0;T^{1/\theta}}\,ds
	\\
	&
	\le C\int_\eta^\infty s^{-\frac{N}{\theta}\left(\frac{1}{r}-\frac{1}{q_1}\right)}\|\mu\|_{r,0;T^{1/\theta}}\,ds
	\le C\eta^{1-\frac{N}{\theta}\left(\frac{1}{r}-\frac{1}{q_1}\right)}\,
	\|\mu\|_{r,0;T^{1/\theta}}.
\end{align*}
These together with \eqref{eq:4.2} and \eqref{eq:4.9} imply that
\begin{equation}
\label{eq:5.3}
\begin{split}
K(f_z,\lambda; L^{q_0},L^{q_1})
 & \le\|f_{z,0}^\eta\|_{L^{q_0}}+\lambda\|f_{z,1}^\eta\|_{L^{q_1}}\\
 & \le C\left(\eta^{1-\frac{N}{\theta}\left(\frac{1}{r}-\frac{1}{q_0}\right)}+\lambda\eta^{1-\frac{N}{\theta}\left(\frac{1}{r}-\frac{1}{q_1}\right)}\right)\|\mu\|_{r,0;T^{1/\theta}}
 \end{split}
\end{equation}
for $z\in{\mathbb R}^N$, $\lambda>0$, and $\eta>0$.
Taking $\eta>0$ so that 
\begin{equation}
\label{eq:5.4}
\eta^{1-\frac{N}{\theta}\left(\frac{1}{r}-\frac{1}{q_0}\right)}=\lambda^\frac{1}{2},
\end{equation}
since $1/(pr)=1/r-\theta/N$, by \eqref{eq:5.1} we have
\begin{equation}
\label{eq:5.5}
	\eta^{1-\frac{N}{\theta}\left(\frac{1}{r}-\frac{1}{q_1}\right)}
	=\eta^{1-\frac{N}{\theta}\left(\frac{1}{r}-\frac{2}{pr}+\frac{1}{q_0}\right)}
	=\eta^{-1+\frac{N}{\theta}\left(\frac{1}{r}-\frac{1}{q_0}\right)}
	=\lambda^{-\frac{1}{2}}.
\end{equation}
Then, by \eqref{eq:5.3}, \eqref{eq:5.4}, and \eqref{eq:5.5}
we see that
$$
K(f_z(t),\lambda; L^{q_0},L^{q_1})
\le C\lambda^{\frac{1}{2}}\|\mu\|_{r,0;T^{1/\theta}}
$$
for $z\in{\mathbb R}^N$, $t\in(0,T)$, and $\lambda>0$. 
By Lemma~\ref{Lemma:4.1} with $\alpha=0$, \eqref{eq:1.7}, \eqref{eq:4.3}, and \eqref{eq:5.1} 
we apply Lemma~\ref{Lemma:4.1} with $\kappa=1/2$ to obtain 
\begin{align*}
\sup_{0<t<T}\bigg\|\int_0^t S(s)\mu\,ds\,\bigg\|'_{pr,0;T^{1/\theta}}
 & \le C\sup_{0<t<T}\sup_{z\in{\mathbb R}^N}\|f_z(t)\|_{L^{pr,\infty}}
 \le C\sup_{0<t<T}\sup_{z\in{\mathbb R}^N}\|f_z(t)\|_{(L^{q_0},L^{q_1})_{1/2,\infty}}\\
 & =C\sup_{0<t<T}\sup_{z\in{\mathbb R}^N}\sup_{\lambda>0}
\left\{\lambda^{-\frac{1}{2}}K(f_z(t),\lambda; L^{q_0},L^{q_1})\right\}\\
 & \le C\|\mu\|_{r,0;T^{1/\theta}}.
\end{align*}
Thus we obtain the desired inequality. 
The proof is complete. 
$\Box$
\vspace{8pt}

Let $\varepsilon>0$ and fix $T\in(0,\infty]$. 
Define
$$
Y_T:=C((0,T): L_{{\rm ul}}^{pr,\infty}).
$$
Set $C^{**}=2C_{**}$, where $C_{**}$ is given in Proposition~\ref{Proposition:5.1}.
For any $u\in Y_T$ and $\varepsilon>0$, we define 
\begin{equation}
\label{eq:5.6}
Y_T(C^{**}\varepsilon):=\left\{u\in Y_T\,:\,\sup_{0<t<T}\,\|u(t)\|'_{pr,0;T^{1/\theta}}\le C^{**}\varepsilon\right\}.
\end{equation}
For any $u$, $v\in Y_T$, set 
$$
d_Y(u,v) :=\sup_{0<t<T}\,\|u(t)-v(t)\|'_{pr,0;T^{1/\theta}}.
$$
Then $(Y_T,d_Y)$ is a complete metric space and $Y_T(C^{**}\varepsilon)$ is closed in $(Y_T,d_Y)$. 
For the proof of Theorem~\ref{Theorem:1.4}, we prepare the following lemma.
\begin{lemma}
\label{Lemma:5.1}
Let $\theta\in(0,2]$ with $N>\theta$, $T\in(0,\infty]$, $p>p_*$, and $r=r_*$. 
Then there exists $C>0$, independent of $T$, such that 
$$
	d_Y({\mathcal F}(u),{\mathcal F}(v))\le C\varepsilon^{p-1}\,d_Y(u,v)
	\quad\mbox{for}\quad u,v\in Y_T(C^{**}\varepsilon),\,\,\varepsilon>0,
$$
where ${\mathcal F}$ is as in \eqref{eq:4.18}.
\end{lemma}
{\bf Proof.}
Let $u$, $v\in Y_T(C^{**}\varepsilon)$. 
Let $0<s<t<T$. 
It follows from Lemmas~\ref{Lemma:2.1}, \ref{Lemma:2.2}, \ref{Lemma:3.5}, and \eqref{eq:4.19} that
$$
	\begin{aligned}
 	& 
	\|F_p(u(s))-F_p(v(s))\|_{r,0;T^{1/\theta}}
	\\
 	& 
	\le C\|w(s)\|_{(pr)/(p-1),0;T^{1/\theta}} \|u(s)-v(s)\|_{pr,0;T^{1/\theta}}
	\\
 	& 
	\le C\left(\|u(s)\|^{p-1}_{pr,0;T^{1/\theta}}+\|v(s)\|^{p-1}_{pr,0;T^{1/\theta}}\right)
	\|u(s)-v(s)\|_{pr,0;T^{1/\theta}},
	\end{aligned}
$$
where $w(x,s):=p(|u(x,s)|^{p-1}+|v(x,s)|^{p-1})$.
By \eqref{eq:5.6} we obtain
\begin{equation}
\label{eq:5.7}
	\sup_{0<s<t}\,\|F_p(u(s))-F_p(v(s))\|_{r,0;T^{1/\theta}}
	\le C(C^{**}\varepsilon)^{p-1}d_Y(u,v).
\end{equation}

For any $z\in{\mathbb R}^N$ and $\eta>0$, 
let $g_z$, $g_{z,0}^\eta$, and $g_{z,1}^\eta$ be as in \eqref{eq:4.21} and \eqref{eq:4.22}. 
Let $q_0$, $q_1$ be as in \eqref{eq:5.1}.
Then, by Proposition~\ref{Proposition:3.2}~(2), \eqref{eq:5.2}, and \eqref{eq:5.7} we see that
$$
	\begin{aligned}
	\|g_{z,0}^\eta(t)\|_{L^{q_0}}
	 & \le
	\int_0^\eta\chi_{(0,t)}(s)\,\bigg\|\chi_{B(z,T^{1/\theta})}S(s)\{F_p(u(t-s))-F_p(v(t-s))\}
	\bigg\|_{L^{q_0}}\,ds
	\\
	&
	\le C\int_0^\eta\chi_{(0,t)}(s)\, s^{-\frac{N}{\theta}\left(\frac{1}{r}-\frac{1}{q_0}\right)}
	\|F_p(u(t-s))-F_p(v(t-s))\|_{r,0;T^{1/\theta}}\,ds
	\\
	&
	\le C\int_0^\eta s^{-\frac{N}{\theta}\left(\frac{1}{r}-\frac{1}{q_0}\right)}\,ds\,
	\cdot\sup_{0<s<t}\,\|F_p(u(s))-F_p(v(s))\|_{r,0;T^{1/\theta}}
	\\
	&
	\le
	C(C^{**}\varepsilon)^{p-1}\eta^{1-\frac{N}{\theta}\left(\frac{1}{r}-\frac{1}{q_0}\right)}d_Y(u,v),\\
	\|g_{z,1}^\eta(t)\|_{L^{q_1}}
	&
	\le
	\int_{\eta}^\infty\chi_{(0,t)}(s)\,\bigg\|\chi_{B(z,T^{1/\theta})}S(s)\{F_p(u(t-s))-F_p(v(t-s))\}
	\bigg\|_{L^{q_1}}\,ds
	\\
	&
	\le C\int_{\eta}^\infty\chi_{(0,t)}(s)\, s^{-\frac{N}{\theta}\left(\frac{1}{r}-\frac{1}{q_1}\right)}
	\|F_p(u(t-s))-F_p(v(t-s))\|_{r,0;T^{1/\theta}}\,ds
	\\
	&
	\le C\int_{\eta}^\infty s^{-\frac{N}{\theta}\left(\frac{1}{r}-\frac{1}{q_1}\right)}\,ds\,
	\cdot\sup_{0<s<t}\,\|F_p(u(s))-F_p(v(s))\|_{r,0;T^{1/\theta}}
	\\
	&
	\le
	C(C^{**}\varepsilon)^{p-1}\eta^{1-\frac{N}{\theta}\left(\frac{1}{r}-\frac{1}{q_1}\right)}d_Y(u,v).
	\end{aligned}
$$
Then, similarly to \eqref{eq:5.3}, we have 
\begin{equation}
\label{eq:5.8}
\begin{split}
K(g_z(t),\lambda; L^{q_0},L^{q_1})
 & \le\|g_{z,0}^\eta(t)\|_{L^{q_0}}+\lambda\|g_{z,1}^\eta(t)\|_{L^{q_1}}\\
 & \le C(C^{**}\varepsilon)^{p-1}\left(\eta^{1-\frac{N}{\theta}\left(1-\frac{1}{q_0}\right)}+\lambda\eta^{1-\frac{N}{\theta}\left(1-\frac{1}{q_1}\right)}\right)d_Y(u,v)
\end{split}
\end{equation}
for $z\in{\mathbb R}^N$, $t\in(0,T)$, $\lambda>0$, and $\eta>0$. 
Taking $\eta>0$ satisfying  \eqref{eq:5.4}, 
by \eqref{eq:5.5} and \eqref{eq:5.8} we obtain 
$$
K(g_z(t),\lambda; L^{q_0},L^{q_1})\le C(C^{**}\varepsilon)^{p-1}\lambda^{\frac{1}{2}}d_Y(u,v)
$$
for $z\in{\mathbb R}^N$, $t\in(0,T)$, and $\lambda>0$. 
Therefore we deduce from Lemma~\ref{Lemma:3.4},
Lemma~\ref{Lemma:4.1} with $\alpha=0$, \eqref{eq:1.8}, \eqref{eq:4.3}, and \eqref{eq:5.1} that 
\begin{align*}
d_Y({\mathcal F}(u),{\mathcal F}(v))
 & =\sup_{t\in(0,T)}\sup_{z\in{\mathbb R}^N}\|g_z(t)\|'_{L^{pr,\infty}(\log L)^0}
\le C\sup_{t\in(0,T)}\sup_{z\in{\mathbb R}^N}\|g_z(t)\|_{(L^{q_0},L^{q_1})_{1/2,\infty}}\\
 & =C\sup_{t\in(0,T)}\sup_{z\in{\mathbb R}^N}\sup_{\lambda>0}\left\{\lambda^{-\frac{1}{2}}K(g_z(t),\lambda; L^{q_0},L^{q_1})\right\}\\
 & \le C(C^{**}\varepsilon)^{p-1}d_Y(u,v).
\end{align*}
Then we obtain the desired inequality. Thus Lemma~\ref{Lemma:4.2} follows.
$\Box$
\vspace{5pt}

\noindent
{\bf Proof of Theorem~\ref{Theorem:1.4}.}
Fix $T\in(0,\infty]$ and $r=r_*$. 
Let $\varepsilon>0$ be small enough.
Let $\mu\in L^{r,\infty}_{{\rm ul}}$ be such that 
$\|\mu\|_{r,0;T^{1/\theta}}\le\varepsilon$. 
By Proposition~\ref{Proposition:5.1}, Lemma~\ref{Lemma:5.1}, \eqref{eq:4.18}, and \eqref{eq:5.6} we have
$$
	\begin{aligned}
	\sup_{t\in(0,T)}\left\|{\mathcal F}(u(t))\right\|'_{pr,0;T^{1/\theta}} 
	& 
	\le \sup_{0<t<T}\bigg\|\int_0^tS(s)\mu\,ds\,\bigg\|'_{pr,0,T^{1/\theta}}
	+{d_Y\bigl({\mathcal F}(u),{\mathcal F}(0)\bigr)}
	\\
	& 
	\le C_{**}\varepsilon+C\varepsilon^{p-1}d_Y(u,0)
	 \le C_{**}\varepsilon+C\varepsilon^{p-1}\cdot C^{**}\varepsilon\le C^{**}\varepsilon
	\end{aligned}
$$
for $u\in Y_T(C^{**}\varepsilon)$. 
This implies that ${\mathcal F}(u)\in Y_T(C^{**}\varepsilon)$ for $u\in Y_T(C^{**}\varepsilon)$.
Furthermore, 
taking small enough $\varepsilon>0$ if necessary, 
by Lemma~\ref{Lemma:5.1} we have 
$$
	d_Y({\mathcal F}(u),{\mathcal F}(v))
	\le C\varepsilon^{p-1}\,d_Y(u,v)\le\frac{1}{2}d_Y(u,v)
$$
for $u$, $v\in Y_T(C^{**}\varepsilon)$. 
We apply the contraction mapping theorem to find a unique $u_*\in Y_T(C^{**}\varepsilon)$ such that 
${\mathcal F}(u_*)=u_*$ in $Y_T(C^{**}\varepsilon)$. 
Then $u_*$ is a solution to problem~\eqref{eq:P} in ${\mathbb R}^n\times[0,T)$ 
in the sense of Definition~\ref{Definition:4.1}, with $u_*$ satisfying \eqref{eq:1.13}. 
Thus Theorem~\ref{Theorem:1.4} follows.
$\Box$\vspace{5pt}
\newline
{\bf Proof of Corollary~\ref{Corollary:1.3}.} 
Let $\mu \in L^{r_*,\infty}_{{\rm ul}}(\log L)^\alpha$, where $\alpha>0$. 
For any $T>0$, set $\omega(T):=|B(0,T^{1/\theta})|$. 
Then it follows from (R2) and (R3) that
$$
\left(\chi_{B(z,T^{1/\theta})}\mu\right)^*(s)\le \chi_{[0,\omega(T))}\left(\frac{s}{2}\right)\mu^*\left(\frac{s}{2}\right),\quad s>0.
$$
Then, by ($\Phi$1) and \eqref{eq:1.5}
we have
\begin{align*}
\|\mu\|_{r_*,0;T^{1/\theta}}
& =\sup_{z\in{\mathbb R}^N}\|\chi_{B(z,T^{1/\theta})}\mu\|_{L^{r_*,\infty}}
=\sup_{z\in{\mathbb R}^N}\sup_{s>0}\left\{s^{\frac{1}{r_*}}(\chi_{B(z,T^{1/\theta})}\mu)^*(s)\right\}\\
 & =\sup_{z\in{\mathbb R}^N}\sup_{0<s<2\omega(T)}\left\{\Phi(s^{-1})^{-\frac{\alpha}{r_*}}\cdot s^{\frac{1}{r_*}}\Phi(s^{-1})^{\frac{\alpha}{r_*}}(\chi_{B(z,T^{1/\theta})}\mu)^*(s)\right\}\\
 & \le \Phi((2\omega(T))^{-1})^{-\frac{\alpha}{r_*}}\|\mu\|_{r_*,\alpha;T^{1/\theta}}\to 0\quad\mbox{as}\quad T\to+0.
\end{align*}
Therefore, by Theorem~\ref{Theorem:1.4} we find a solution $u$ to problem~\eqref{eq:P} in ${\mathbb R}^N\times(0,T)$ 
for small enough $T>0$.  
Thus Corollary~\ref{Corollary:1.3} follows.
$\Box$\vspace{5pt}
\newline
{\bf Proof of Corollary~\ref{Corollary:1.4}.} 
It follows from \eqref{eq:1.14} that $\mu_c\in L^{r_*,\infty}$. 
Then Corollary~\ref{Corollary:1.4} follows from Theorem~\ref{Theorem:1.4}.
$\Box$ 
\vspace{5pt}

\noindent
{\bf Acknowledgements.} 
K. I. and T. K. were supported in part by JSPS KAKENHI Grant Number JP19H05599.
T. K. was supported in part by JSPS KAKENHI Grant Number JP22KK0035.
R. T. was supported in part by JSPS KAKENHI Grant Numbers JP22K03388, JP21H00993, JP20H01814, JP18KK0072.
\vspace{5pt}

\noindent
{\bf  Conflict of Interest.}
The authors state no conflict of interest. 

\end{document}